\documentclass[11pt,a4paper,reqno]{amsart}
\usepackage{a4wide}
\usepackage{amsfonts,amsmath,amssymb,amsthm,enumerate,color,hyperref,bbm,amsthm,epsfig,amscd,graphicx}
\usepackage[latin1]{inputenc}

\newcommand{\eps}{\varepsilon}
\newcommand{\id}{{\rm id}}
\newcommand{\dd}{\;{\rm d}}
\newcommand{\Fhls}{{\mathcal F}_{\rm PKS}}
\newcommand{\Fcfd}{{\mathcal H}_{\lambda}}
\newcommand{\F}{{\mathcal F}_{\rm PKS}}
\newcommand{\R}{{\mathbb R}}
\newcommand{\RR}{{\mathbb R}}
\newcommand{\GG}{{\mathcal G}}
\newcommand{\N}{{\mathbb N}}
\newcommand{\NN}{{\mathbb N}}
\newcommand{\LL}{{\rm L}}
\newcommand{\PP}{{\mathcal P}}
\newcommand{\W}{{\mathcal W}}
\newcommand{\rhohls}{\bar \varrho_\lambda}
\newcommand{\Fcfdr}{{\mathcal H}_{\lambda,\delta}}

\newcommand{\EE}{{\mathcal K}}

\renewcommand{\AA}{\mathcal A}

\newcommand{\quid}[1]{\noindent{\bf Open question: }\emph {#1}}
\newtheorem{theorem}{Theorem}
\newtheorem{proposition}{Proposition}

\theoremstyle{remark}
\newtheorem{remark}{Remark}

\title[On the parabolic-elliptic Patlak-Keller-Segel system]{On the parabolic-elliptic Patlak-Keller-Segel system in dimension 2 and higher}
\author{Adrien Blanchet}
\address{Toulouse School of Economics (GREMAQ, Universit\'e de Toulouse), 21 All\'ee de Brienne, F-31000 Toulouse, France}
\email{Adrien.Blanchet@univ-tlse1.fr}
\date{\today}
\thanks{Work supported by the project EVaMEF ANR-09-JCJC-0096-01}
\begin{document}
\begin{abstract}
This review is dedicated to recent results on the 2d parabolic-elliptic Patlak-Keller-Segel model, and on its variant in higher dimensions where the diffusion is of critical porous medium type. Both of these models have a critical mass $M_c$ such that the solutions exist globally in time if the mass is less than $M_c$ and above which there are solutions which blowup in finite time. The main tools, in particular the free energy, and the idea of the methods are set out. A number of open questions are also stated.
\end{abstract}
\maketitle
\section{Biological background}
Chemo-taxis is defined as a move of an organism along a chemical concentration gradient. Bacteria can produce this chemo-attractant themselves, creating thus a long-range non-local interaction between them. We are interested in a very simplified model of aggregation at the scale of cells by chemo-taxis: myxamoebaes or bacterias experience a random walk to spread in the space and find food. But in starvation conditions, the dictyostelium discoideum emit a chemical signal: cyclic adenosine monophosphate (cAMP). They move towards a higher concentration of cAMP. Their behaviour is thus the result of a competition between a random walk-based diffusion process and a chemo-taxis-based attraction. It was noticed experimentally that if there are enough bacteria they aggregate whereas if they are not enough they go on spreading in a chemo-tactically inert environment, {\it e.g.}~\cite{brenner}. The typical time scale for the spreading of bacteria on a petri dish is around one day, and a few minutes for the concentration. This concentration phenomenon is the first step for uni-cellular organisms to come together with others and form a multi-cellular organism. It can be seen as a hint on how, during the evolution of species, the passage from uni-cellular organisms to more complex structure was achieved. It is also a paradigm model for pattern formation of cells for meiose, embryo-genesis or angio-genesis.

In nature the dictyostelium discoideum spread on the soil and then come together to form a motile pseudoplasmodium. This slug creeps to a few centimeters below the soil surface where it forms a fruiting body with spores and a stalk. The spores are then blown away by the wind to colonise a new place. Around 20\% of the cells which are in the stalk altruistically sacrifice themselves to allow the species to survive. They are an excellent example of social behavior with outstanding coordination and sense of sacrifice for the benefit of the species.
\section{The 2d parabolic-elliptic Patlak-Keller-Segel system}
\subsection{The model}
The first mathematical attempt to model this aggregation phenomenon is often granted to Evelyn Fox Keller and Lee A. Segel in~\cite{Keller-Segel-70} but the same model was earlier described by Patlak in~\cite{Pat53}. We consider the following simplified version was given by~\cite{JL92}:
\begin{equation}\label{eq:KS}
\left\{
\begin{array}{ll}
 \displaystyle \frac{\partial \rho}{\partial t}=\Delta \rho-\nabla \cdot (\rho\nabla c)
 &\quad  x\in\mathbb{R}^2,\;t>0\,,\vspace{.3cm}\\
-\Delta c=\rho&\quad  x\in\mathbb{R}^2,\;t>0\,,\vspace{.3cm}\\
\rho(\cdot,t=0)=\rho_0\ge 0&\quad x\in\mathbb{R}^2\,.
\end{array}
\right.
\end{equation}
Here $\rho$ represents the cell density and $c$ the concentration of chemo-attractant.

As the solution to the Poisson equation $-\Delta c=\rho$ is given up to a harmonic function, we choose the one given by $c=G*\rho$ where 
\begin{equation*}
G(|x|):=-\frac1{2\pi}\log|x|\;.
\end{equation*}
The Patlak-Keller-Segel system~\eqref{eq:KS} can thus be written as a non-local parabolic equation:
\begin{equation*}
  \frac{\partial \rho}{\partial t}=\Delta \rho - {\rm div}(\rho\nabla G*\rho)\quad  \mbox{in }\,(0,+\infty) \times \mathbb R^2\;.
\end{equation*}
Also note that the mass is conserved
\begin{equation*}
\int_{\mathbb R^2}\rho(t,x)\dd x=\int_{\mathbb R^2}\rho_0(x)\dd x=:M\;.
\end{equation*}

S. Childress and J. Percus and V. Nanjundiah conjectured in~\cite{MR632161,nanjundiah} that this system displays the existence of a critical mass above which the cells aggregate and below which they do not.,For a complete review of the early literature, the interested reader could beneficially consult~\cite{H03,H03a}. For a more recent references see~\cite{perthame} or~\cite{HP}.\medskip

Except when it is clearly indicated, in this article we will assume that the initial condition satisfies:
\begin{equation*}
(1+|x|^2)\,\rho_0\in \LL^1_+(\RR^2)\quad\mbox{and}\quad \rho_0\,\log \rho_0\in
\LL^1(\RR^2)\;.
\end{equation*}
%%%%%%%%%%%%%%%%%%%%%%
\subsection{Blowup}
Consider a smooth solution to the Patlak-Keller-Segel system~\eqref{eq:KS} we formally compute a virial identity
\begin{equation}\label{eq:2moment}
  \frac \dd{\dd t}\int_{\mathbb{R}^2}|x|^2 \rho(t,x)\dd x = 4M\big(
1-\frac{M}{8\pi} \big)\,.
\end{equation}
So that all the solutions with mass bigger than $8\pi$ and finite 2-moment cannot be global in time:
%---------------------------------------------------
\begin{theorem}[Blowup]
Let $\rho$ be a solution to the Patlak-Keller-Segel system~\eqref{eq:KS} and $[0, T^*)$ its maximal
interval of existence. If $M>8\pi$,
then
\[
T^*\le \frac{2\pi}{M(M-8\pi)}\int_{\mathbb{R}^2} |x|^2  \rho_0(x)\dd x\,,
\] 
and $\rho(t,\cdot)$ converges, up to extraction of a sub-sequence, as $t\to T^*$ to a measure which is not in $L^1(\mathbb{R}^2)$.
\end{theorem}
%---------------------------------------------------
This virial method was already used in the dispersive equations community. However, this kind of method is non-constructive and gives no hint on the nature of the blowup. Another non-constructive proof of this non global-in-time existence will be given in Remark~\ref{benasque} and will give similar results even in the case when the solutions are of 2-moment infinite.

The proof that the blowup occurs in an isolated point was initiated in~\cite{nagaisenbasuzuki,medina}. When the solution is radially symmetric, the blowup mechanism is not clear and the different conjectures are as follows:
There exists a solution to the Patlak-Keller-Segel system~\eqref{eq:KS} that yields a concentration of an amount of mass $8\pi$ at the origin. Moreover, the asymptotics of this solution near the origin is given by
\begin{equation}\label{mechablo}
  \rho(t,x) \sim\frac{1}{L(t)^2}\frac{8}{\left( 1 + |y|^2\right)^2}\quad\mbox{and}\quad c(t,x)=-2\log(1+y^2)
\end{equation}
for $t \to T^-$ where $y=r/L(t)$ and as $t$ goes to $\infty$. The different approach lead to a different expression for $L$. The most recent approach is due to~\cite{Luxx} where the author makes strong analogies between the Patlak-Keller-Segel system and the non-linear Sch\"odinger equation. We will come back to these analogies in the concluding remarks in Section~\ref{sec:conclusion}. P. Lushnikov uses an appropriate gauge transformation to a new time dependent variable around the self-similar solution to obtain a linearisation operator of self-adjoint form. The study of the spectrum of the linearised operator allows him to derive a set of amplitude equations for the coefficients of this expansion. He thus obtain
\begin{equation}\label{defl}
  L(t)=2\sqrt{T-t}\,e^{-1-\gamma/2} \,e^{\sqrt{-\log(T-t)/2}}\left[1+O\left( \frac{\left[\log\left(-\log(T-t) \right) \right]^3}{\sqrt{-2\log(T-t)}}\right) \right]
\end{equation}
$\gamma$ being the Euler constant. In~\cite{velazquez}, a similar result was obtained using matched asymptotic expansion techniques  but with exponent $1$ in the log-log term in~\eqref{defl}. Whereas~\cite{MR1415081} conjectures
\begin{equation*}
  L(t) \sim \sqrt{T-t}\,e^{\sqrt{-\log(T-t)/2}} \left(-\log(T-t) \right)^{(-\log(T-t))^{-1/2}/4}\;.
\end{equation*}
If the numerics performed in~\cite{SC02} agree on the leading order term $\sqrt{T-t}$, C.~Sire and P.-H.~Chavanis predict that
\begin{equation*}
  L(t) \sim \sqrt{T-t}\,e^{-\sqrt{-\log(T-t)\log(-\log(T-t))/2}/2}\;.
\end{equation*}

\quid{All those estimates are formal. Justification of this result is very challenging, as is the numerical track of the next order term in $L$.}
\medskip

Concerning the continuation of the solution after blowup, the usual idea is to define a sequence of approximate problems containing a small parameter $\eps>0$ which have global-in-time solutions and approach the original problem when $\eps$ goes to $0$. The behaviour of the solution to the approximate model has to be close to the one of the Patlak-Keller-Segel system~\eqref{eq:KS}, except close to the singularities. The blowup indicates that the approximate problem is not valid anymore close to the singularity. In~\cite{MR2068668,MR2068667}, using matched asymptotic expansions J. Vel{\'a}zquez describes in a rather detailed manner the formation and motion of some regions where the mass concentrates. He also proves local existence of the solution as long as there is no formation of another Dirac mass and no collision of Dirac masses occur. The approximate problem is
\begin{equation*}
\frac{\partial \rho}{\partial t}=\Delta \rho-\nabla \cdot (\Phi_\eps(\rho)\nabla G*\rho)
 \quad  x\in\mathbb{R}^2,\;t>0\,,
\end{equation*}
where $\Phi_\eps(\rho):=\eps^{-1}\Phi(\eps\rho)$ where $\eps$ is a small parameter and $\Phi$ is an increasing bounded function satisfying
\begin{equation*}
  \Phi(s)=s-\alpha s^2 + \cdots \;\mbox{as $s \to 0$ and }\Phi(s)\sim L\;\mbox{as $s \to \infty$}\;,
\end{equation*}
where $L>0$ is a given number. This model prevents overcrowding as the chemo-tactic function $\Phi_\eps(\rho)$ saturates at the constant value $L$.

In~\cite{DS}, J.~Dolbeault and C.~Schmeiser define measure valued densities to give a sense to generalized global-in-time solutions for any mass. This extends the solution concept after blow-up. They also show that the choice of a solution concept after blow-up is not unique and depends on the type of regularization. Their regularisation is different from the one chosen above. The regularised problem they consider is
\begin{equation}\label{eq:ksepsds}
\frac{\partial \rho}{\partial t}=\Delta \rho-\nabla \cdot (\rho\nabla G_\eps*\rho)
 \quad  x\in\mathbb{R}^2,\;t>0\,,
\end{equation}
 where $G_\eps(x):=-\log(|x|+\eps)/2\pi$.
%---------------------------------------------------
 \begin{theorem}[Generalised solution,~\cite{DS}]
   For every $T > 0$, as $\eps \to 0$, a sub-sequence of solutions $\eps$ to~\eqref{eq:ksepsds} converges tightly and uniformly in time to a time dependent measure $\rho(t)$.
There exists $\nu(t)$ such that $(\nu,\rho) $ is a generalized solution in the distributional sense of
\begin{equation}\label{eq:geneks}
  \frac{\partial \rho}{\partial t}+{\rm div}(j[\rho,\nu] - \nabla \rho)=0
\end{equation}
where the convective flux $j[\rho,\nu]$ is supported in the support of $\nu$ and is given by
\begin{multline*}
  \int_0^T \int_{\R^2}\varphi(t,x)j[\rho,\nu](t,x)\dd x \dd t= - \frac1{4\pi}\int_0^T \int_{\R^2}\nu(t,x)\nabla \varphi(t,x)\dd x \dd t\\ -\frac1{4\pi}\int_0^T \int_{\R^4} \left(\varphi(t,x)- \varphi(t,y)\right)K(x-y)\rho(t,x)\rho(t,y)\dd x \dd y \dd t
\end{multline*}
for any $\varphi \in \mathcal C^1_b((0,T) \times \R^2)$ with
\begin{equation*}
\left\{
  \begin{array}{ll}
    \displaystyle \frac{x}{|x|^2}&\quad \mbox{for $x \neq 0$}\vspace{.1cm}\\
    0&\quad \mbox{for $x = 0$}\;.
  \end{array}
\right.
\end{equation*} 
 \end{theorem}
%--------------------------------------------------
If $\rho$ does not charge points then the additional default measure $\nu$ vanishes and $j[\rho,0]=\rho \nabla G*\rho$, so that~\eqref{eq:geneks} is  generalisation of the Patlak-Keller-Segel system~\eqref{eq:KS}. They also obtain a strong formulation when the generalised solution is assumed to be the sum of a regular part and of Dirac masses:
%---------------------------------------------------
 \begin{theorem}[Strong formulation,~\cite{DS}]
Assume that the generalised solution to~\eqref{eq:geneks} has the form
\begin{equation}\label{assuinde}
    \rho(t,x) = \sum_{j \in N}M_j(t) \delta(x-x_j(t)) + \rho_{\rm reg}(t,x)\;.
  \end{equation}
Then as long as there is no formation of a new Dirac mass and that there is no collision of Dirac masses
\begin{equation*}
  \frac{\partial \rho_{\rm reg}}{\partial t}= \Delta \rho_{\rm reg}- \nabla\left( \rho_{\rm reg} \nabla G*\rho_{\rm reg}-\frac{1}{2\pi} \nabla\rho_{\rm reg}\sum_{j _in N}M_j(t) \frac{x-x_j(t)}{|x-x_j(t)|^2}\nabla\rho_{\rm reg} \right)
\end{equation*}
with
\begin{equation*}
  \dot M_i(t) = M_i(t)\rho_{\rm reg}(t,x_i(t))
\end{equation*}
and
\begin{equation*}
  \dot x_i(t)= \nabla G*\rho_{\rm reg}(t,x_i(t)) -\frac{1}{2\pi} \sum_{j\in N,i\neq j}M_j(t) \frac{x-x_j(t)}{|x-x_j(t)|^2}\;.
\end{equation*}
 \end{theorem}
%--------------------------------------------------
A similar result was formally obtained in~\cite{MR2068668}, with the last equation replaced by 
\begin{equation*}
  \dot x_i(t)=\Gamma(M_i(t)) \left(-\frac{1}{2\pi} \sum_{j\in N,i\neq j}M_j(t) \frac{x-x_j(t)}{|x-x_j(t)|^2} + \nabla G*\rho_{\rm reg}(t,x_i(t)) \right)
\end{equation*}
where $\Gamma(M)$ is a mean value of the derivative $\Phi$ is such that $0<\Gamma(M)<1$, $\Gamma(8\pi)=1$ and $\Gamma(\infty)=0$ and can be explicitly described, see~\cite[Equation (3.45)]{MR2068668}. Note that this means that the limit depends on the choice of the regularisation.

In~\cite{ChMannella}, for the case of two particles, the authors prove that when $t$ is large enough the solution is made of a Dirac peak of mass $M_0(t)$ surrounded by a dilute halo containing the remaining mass whose dynamical evolution is described by a Fokker-Planck equation. Therefore, they neglect the self-gravity of the halo and prove that the mass of the Dirac peak saturates to $M$ algebraically rapidly as 
\begin{equation*}
  1-\frac{M_0(t)}{M}\sim t^{-a}\;.
\end{equation*}

\quid{Actually the assumption~\eqref{assuinde} is valid only between two blowup events or between two collisions of Dirac masses. We expect that in the end the solution is made of one Dirac mass with all the mass but we are still missing such a rigorous theory.}
%%%%%%%%%%%%%%%%%%%
\subsection{Global existence}\label{subcritical}
%-------------
\subsubsection{{\it A priori} estimates}
The first important answer to Nanjundiah's conjecture on global existence {\it vs} blowup was given by W.~J\"ager and S.~Luckhaus in~\cite{JL92}: the natural idea is to regularise the Green kernel and to pass to the limit. The solutions have mass $M$ so that the loss of compactness can come either from concentration or vanishing. By~\eqref{eq:2moment}, the second moment remains bounded so that the main problem is to control the concentration of mass. W.~J\"ager and S.~Luckhaus tried to obtain a bound on the entropy $\int_{\mathbb{R}^2}\rho\log \rho\dd x$ by differentiating it and using an integration by parts and the equation for $c$, we obtain:
\begin{equation*}
\frac \dd{\dd t}\int_{\mathbb{R}^2}\rho\log \rho\dd x=-4\int_{\mathbb{R}^2}\left|\nabla\sqrt \rho\right|^2 \dd x
  + \int_{\mathbb{R}^2} \rho^2 \dd x \,.
\end{equation*}
Applying the Gagliardo-Nirenberg-Sobolev inequality:
\begin{equation}\label{GNS} \int_{\R^2}|u|^4\dd x\le C_{\rm GNS} \int_{\R^2}|\nabla u|^2 \dd x\int_{\R^2}|u|^2\dd x\quad\forall\; u\in H^1(\mathbb{R}^2)\,.
\end{equation}
to $u=\sqrt \rho$, we have
\begin{equation*}
\frac \dd{\dd t}\int_{\mathbb{R}^2}\rho\log \rho\dd x \le  \left[-4+MC_{\rm GNS} \right]\int_{\mathbb{R}^2}\left|\nabla\sqrt \rho\right|^2 \dd x 
\end{equation*}
So that the entropy is non-increasing if $M\le 4/C_{\rm GNS}\approx 1.862... \times 4\pi<8\pi$. They hence obtained global existence in this case together with important propagation of the $L^p$-estimates.\medskip

We can indeed improve this result by using the following free energy:
\[
\F[\rho]:=\int_{\mathbb{R}^2}\rho\log \rho\dd x - \frac{1}2\int_{\mathbb{R}^2}\rho c\dd x\,.
\]
A simple formal calculation shows that for  all
$u \in C_c^\infty(\R^2)$ with zero mean,
$$
\lim_{\epsilon\to 0}\frac{1}{\epsilon}\left(\Fhls[\rho+\epsilon u]
-\Fhls[\rho]\right) = \int_{\R^2}\frac{\delta
{\Fhls}[\rho]}{\delta  \rho}(x) \, u(x)\dd x
$$
where
\begin{equation*}
\frac{\delta {\Fhls}[\rho]}{\delta  \rho}(x) := \log \rho(x)-G*\rho(x) \ .
\end{equation*}
It is then easy to see that the Patlak-Keller-Segel system~\eqref{eq:KS} can be rewritten as
\begin{equation}\label{wasgrad}
\frac{\partial \rho}{\partial t}(t,x)  = {\rm
div}\left(\rho(t,x)\nabla \left[ \frac{\delta
{\Fhls}[\rho(t)]}{\delta  \rho}(x) \right]\right)\ .
\end{equation}
It follows that at least along well-behaved solutions to the Patlak-Keller-Segel system~\eqref{eq:KS},
\begin{equation*}
\frac{{\dd}}{{\dd}t}{\Fhls}[\rho(t)] = - \int_{\R^2}
\rho(t,x)\left|\nabla  \left[\frac{\delta {\Fhls}[\rho(t)]}{\delta
\rho}(x)\right]\right|^2\dd x \ .
\end{equation*}
Or equivalently
\begin{equation*}
\frac \dd{\dd t}\F[\rho(t,\cdot)]=-\int_{\mathbb{R}^2}\rho\left|\nabla\left(\log \rho - c\right)\right|^2\dd x \,.
\end{equation*}
In particular, along such solutions, $t \mapsto {\Fhls}[\rho(t)]$
is monotone non--increasing.

 The gap between the $4/C_{GNS}$ and $8\pi$ was not filled before~\cite{DP} when was made the link with the logarithmic Hardy-Littlewood-Sobolev: Let $f$ be a non-negative
function in $L^1(\mathbb{R}^2)$ such that $f\log f$ and $f\log
(1+|x|^2)$ belong to $L^1(\mathbb{R}^2)$. If $\int_{\mathbb{R}^2}f
\dd x=M$, then
\begin{equation}\label{eq:HLS}
\int_{\mathbb{R}^2}f\log f\,\dd x+\frac{2}{M}
\iint_{\mathbb{R}^2\times\mathbb{R}^2}f(x)f(y) \log|x-y|\dd x
\dd y\geq -\;C(M)\,,
\end{equation}
with $C(M):=M(1+\log\pi-\log M)$. Moreover the minimisers of~\eqref{eq:HLS} are the translations of
 \begin{equation*}
  {\bar\rho_{\lambda}(x):=\frac{M}{\pi}\frac{\lambda}{\left(\lambda + |x|^2\right)^{2}}\;.} 
 \end{equation*}\medskip

Using the monotonicity of $\Fhls[\rho]$ and the Logarithmic Hardy-Littlewood-Sobolev inequality it is easy to see that
\begin{align}\label{bdp}
\Fhls[\rho] &= \frac{M}{8\pi}\left( \int_{\R^2}\rho(x)\,\log
\rho(x)\dd x+\frac{2}{M} \iint_{\R^2\times\R^2}\rho(x) \log|x-y| \rho(y)\dd x \dd y\right)\nonumber\\
 &\quad+ \left(1 -   \frac{M}{8\pi}\right)\int_{\R^2}\rho(x)\,\log \rho(x)\dd x\nonumber\\
 &\ge - \frac{M}{8\pi}\,C(M) + \left(1 -   \frac{M}{8\pi}\right)\int_{\R^2}\rho(x)\,\log \rho(x)\dd x\ .
\end{align}
It follows that for solutions $\rho$ of the Patlak-Keller-Segel system~\eqref{eq:KS},
\begin{equation}
  \label{eq:entro}
  \int_{\R^2}\rho(t,x)\log \rho(t,x)\dd x \le
\frac{8\pi {\Fhls}[\rho_0]  - M\,C(M)}{8\pi - M}\ .
\end{equation}
Therefore, for $M< 8\pi$, the {\rm entropy} stays bounded uniformly in time. This precludes the collapse of mass into a point mass for such initial data.\medskip

Coming back to the super-critical mass case, it is worth noticing that for a given $\rho$, if we set $\rho_\lambda(x)=\lambda^{-2}\rho(\lambda^{-1}x)$ then
\begin{equation}
  \label{eq:scalingf}
  \F[\rho_\lambda]=\F[\rho]-2M \left(1-\frac{M}{8\pi}\right)
\log\lambda \,.
\end{equation}
So that as a function of $\lambda$, $\F[\rho_\lambda]$ is bounded from below if $M < 8\pi$, and not bounded from below if $M > 8\pi$.
\medskip

As an alternative to the regularisation/passing to the limit procedure, another conceited but smart way to prove the global existence is to use the gradient flow interpretation in the Wasserstein metric. For this purpose we need to introduce a few elements of optimal transport. We denote $\PP(X)$ the set of probability measure on $X$ and $\PP_2(X)$ the subset of probability measures with finite second moments. We say that $T$ {\it transports} $\mu$ onto $\nu$, and denote $T\#\mu=\nu$, if
\begin{equation*}
\int_{\R^2} \zeta(T(x)) \dd \mu(x) = \int_{\R^2} \zeta(y)
\dd \nu(y) \quad \forall \zeta \in \mathcal C_b^0(\R^2)\;.
\end{equation*}
We can define the 2-Wasserstein distance,  $\W_2$ in $\PP(X) \times \PP(Y)$ by
\begin{equation*}
\W_2(\mu, \nu)^2 =\inf_{\Pi\in\Gamma} \iint_{X \times
Y} |x-y|^2 \dd \Pi(x, y) \ ,
\end{equation*}
where $\Pi$ runs over the set of joint probability measures in $X\times Y$ with first marginal $\mu$ and second $\nu$.

Using this metric, we can see the Patlak-Keller-Segel system~\eqref{eq:KS} as a gradient flow of the free energy in the Wasserstein metric:
\begin{equation*}
  \rho_t=-\nabla_{W} \Fhls[\rho(t)]\;.
\end{equation*}
In the sense that we can construct a solution using the minimising scheme, often known as the Jordan-Kinderlehrer-Otto (JKO) scheme: given a time step $\tau$, we define the solution by:
\begin{equation*}
  \rho^{k+1}_\tau \in {\rm argmin}_{\rho \in \mathcal K} \left[\frac{\W^2_2(\rho,\rho^k_\tau)}{2\tau} + \Fhls[\rho] \right]\;,
\end{equation*}
where ${\mathcal K}:= \{\rho\,:\, \int_{\R^2} \rho=M,\; \int_{\R^2} \rho(x) \log \rho(x)\dd x < \infty\; \mbox{and}\; \int_{\R^2}|x|^2 \rho(x)\dd x< \infty\}$.

For the analogy, note that if the metric was Euclidean, the Euler-Lagrange equation associated to
\begin{equation}\label{JKOscheme}
  \rho^{k+1}_\tau \in {\rm argmin} \left[\frac{|\rho-\rho^k_\tau|^2}{2\tau} + \Fhls[\rho] \right]\;,
\end{equation}
would be
\begin{equation*}
  \frac{\rho^{k+1}_\tau-\rho^k_\tau}{\tau} + \Fhls[\rho^{k+1}_\tau] =0\;,
\end{equation*}
which is nothing but the implicit Euler scheme associated to 
\begin{equation*}
  \rho_t=-\nabla \Fhls[\rho(t)]\;.
\end{equation*}

At this point it is convenient to emphasise that the functional $\Fhls$ is not convex, so even the existence of a minimiser is not clear. When the functional is convex, or even displacement convex, general results from~\cite{villani,AmbrosioGigliSavare02p} can be applied. In concrete terms, a functional ${\mathcal G}$ is said to be {\em displacement convex} when the following is true: for any two densities $\rho_0$ and $\rho_1$ of the same mass $M$, let $ \varphi$ be such that $\nabla \varphi\#
\rho_0 = \rho_1$. For $0 < t < 1$, define
\begin{equation*}
\varphi_t(x) = (1-t)\frac{|x|^2}{2} + t\varphi(x)\qquad{\rm and}\qquad \rho_t = \nabla \varphi_t \# \rho_0\ .
\end{equation*}
The {\em displacement interpolation} between $\rho_0$ and $\rho_1$
is the path of densities $t\mapsto \rho_t$, $0\le t \le 1$. Let $\gamma$ be any real number. To say that ${\mathcal G}$ is {\em $\gamma$-displacement convex} means that for all such mass densities $\rho_0$ and $\rho_1$, and all $0\le t \le 1$,
\begin{equation*}
(1-t){\mathcal G}(\rho_0) + t{\mathcal G}(\rho_1) - {\mathcal G}(\rho_t)\ge
\gamma t(1-t)\W_2^2(\rho_0,\rho_1)\ .
\end{equation*}
${\mathcal G}$ is simply {\em  displacement convex} if this is true for $\gamma=0$, and ${\mathcal G}$ is {\em uniformly displacement convex} if this is true for some $\gamma>0$.

We interpolate between the terms of the sequence $\{\rho_\tau^k\}_{k\in\NN}$ to produce a function from $[0,\infty)$ to $L^1(\R^2)$: For each positive integer $k$, let $\nabla\varphi^k$ be the optimal transportation plan with $\nabla\varphi^k \#\rho_\tau^k =  \rho_\tau^{k-1}$.  Then for $(k-1)\tau \leq t \leq k\tau$ we define
$$
\rho_\tau(t)=\left( \frac{t-(k-1)\tau}{\tau} \id +\frac{k\tau
-t}{\tau}\nabla \varphi^k\right) \#\rho^{k}_\tau \ .
$$
%---------------------------------------------------
\begin{theorem}[Convergence of the scheme as $\tau\to0$,~\cite{BCC}]
If $M<8\pi$ then
the family $(\rho_\tau)_{\tau > 0}$ admits a
sub-sequence converging weakly in $L^1(\RR^2)$ to a weak solution to the Patlak-Keller-Segel system.
\end{theorem}
%---------------------------------------------------
In this proof the perturbation of the minimiser has to be done in the Wasserstein setting: Let $\zeta$ be a smooth vector field with compact support, we introduce $T_\eps:={\rm id} + \eps\zeta$. We define
$\overline{\rho_\eps}$ the push-forward perturbation of
$\rho^{n+1}_\tau$ by $T_\eps$:
\[
 \overline{\rho_\eps} = T_\eps \# \rho_\tau^{n+1}\;.
\]
Standard computations give
\begin{multline*}
\int_{\RR^2}  \zeta(x)\frac{x-  \nabla \varphi^{n}(x)}\tau\rho_\tau^{n+1}(x)\dd x
\\= \int_{\RR^2} \left[{\rm div}\zeta(x) -\frac{1}{4\,\pi}\int_{\RR^2}\frac{[\zeta(x)-\zeta(y)]\cdot(x-y)}{|x-y|^2} \rho_\tau^{n+1}(y) \dd y\right] \rho_\tau^{n+1}(x) \dd x\ ,
\end{multline*}
which is the weak form of the Euler-Lagrange equation:
\begin{equation}\label{eq:EL}
  \frac{{\rm id} - \nabla \varphi^{n}}{\tau}\rho_\tau^{n+1} = - \nabla \rho_\tau^{n+1} + \rho_\tau^{n+1} \nabla c_\tau^{n+1}\;.
\end{equation}
Using the Taylor's expansion $\zeta (x)-\zeta\left[\nabla
\varphi^{n}(x)\right] = \left[x-\nabla \varphi^{n}(x)\right]\cdot
\nabla \zeta(x) + O\left[ \left|x-\nabla \varphi^{n}(x)
\right|^2\right]$, we obtain for all$t_2 > t_1 \ge 0$,
\begin{multline}\label{eq:avectau}
\int_{\RR^2} \zeta (x) \left[ \rho_\tau(t_2,x) -
\rho_\tau(t_1,x)\right]\dd x =\int_{t_1}^{t_2}\int_{\RR^2}
\Delta\zeta(x)\,\rho_\tau(s,x)\dd x\dd s + O(\tau^{1/2})\\-\frac1{4\pi}
\int_{t_1}^{t_2}\iint_{\R^2 \times \RR^2}
\rho_\tau(s,x)\,\rho_\tau(s,y)\,\frac{(x-y) \cdot
\left(\nabla\zeta(x)-\nabla\zeta(y)\right)}{|x-y|^2}\dd y\dd
x\ .
\end{multline}
To pass to the limit, the scheme provides some {\it a priori} bounds: Taking $\rho^{n+1}_\tau$ as a test function in~\eqref{JKOscheme} we have:
\begin{equation}
  \label{eq:estime}
  \F[\rho^{n+1}_\tau] + \frac{1}{2\,\tau} \W_2^2(\rho^n_\tau,\rho^{n+1}_\tau) \le \F[\rho^n_\tau]\;.
\end{equation}
As a consequence we obtain an \textit{energy estimate}
\begin{equation*}
 \sup_{n \in \NN}\F[\rho^n_\tau]  \le \F[\rho^0_\tau],
\end{equation*}
which together with~\eqref{bdp} forbids the concentration, and a \textit{total square estimate}
\begin{equation*}
 \frac{1}{2\,\tau}  \sum_{n \in \NN} \W^2_2(\rho^n_\tau,\rho^{n+1}_\tau) \le \F[\rho^0_\tau] - \inf_{n \in \NN}\F[\rho^n_\tau]\;,
\end{equation*}
which rules out vanishing. These two estimates allow to pass to the limit in $\tau$ in~\eqref{eq:avectau}, to obtain:
\begin{multline*}
\int_{\RR^2} \zeta (x) \left[ \rho(t_2,x) -
\rho(t_1,x)\right]\dd x =\int_{t_1}^{t_2}\int_{\RR^2}
\Delta\zeta(x)\,\rho(s,x)\dd x\dd s \\-\frac1{4\pi}
\int_{t_1}^{t_2}\iint_{\R^2 \times \RR^2}
\rho_\tau(s,x)\,\rho(s,y)\,\frac{(x-y) \cdot
\left(\nabla\zeta(x)-\nabla\zeta(y)\right)}{|x-y|^2}\dd y\dd
x\ .
\end{multline*}
Which is the definition of a weak solution.\medskip

We can actually obtain:
%------------------------------
\begin{theorem}[Existence of solution is the subcritical case,~\cite{BDP}] \label{thm:Existence}
If $M < {8\pi}$,
then the Patlak-Keller-Segel system~\eqref{eq:KS} has a global
weak non-negative solution $\rho$ with initial data $\rho_0$ such that
\begin{gather*}
(1+|x|^2+|\log \rho|)  \rho\in L^\infty_{\rm loc}(\mathbb{R}^+,L^1(\mathbb{R}^2))\,,\\
\int_0^t \int_{\mathbb{R}^2} \rho |\nabla \log \rho - \nabla c |^2
\dd x \dd t< \infty\,, \\
\int_{\mathbb{R}^2}|x|^2 \rho(t,x)\dd x = \int_{\mathbb{R}^2}|x|^2
\rho_0(x)\dd x + 4M\big( 1-\frac{M}{8\pi} \big)t
\end{gather*}
for $t>0$. Moreover $\rho \in L^\infty_{\rm loc}((\varepsilon,\infty),L^p(\mathbb{R}^2))$ for any $p\in (1,\infty)$
and any $\varepsilon>0$, and the following inequality holds for
any $t>0$:
\begin{equation*}
\F[\rho(\cdot,t)]+\int_0^t\int_{\mathbb{R}^2}\rho\left|\nabla\left(\log
\rho - c\right)\right|^2\dd x\dd s\le \F[\rho_0]\,.
\end{equation*}
\end{theorem}
%------------------------------
Similar results were first proved in~\cite{nagai} for radially symmetric solutions in a bounded domain with Neumann boundary conditions.\medskip

This notion of free energy solution allows to study the large time behaviour, intermediate asymptotics and convergence to asymptotically self-similar profiles: let $(u_\infty, v_\infty)$ be the unique solution to the Gelfand equation
\begin{equation} \label{eq:stasol}
u_\infty=M \frac{e^{v_\infty-|x|^2/2}}{\int_{\mathbb{R}^2}e^{v_\infty-|x|^2/2}\dd x}= -\Delta v_\infty\,,\quad\mbox{with}\quad
v_\infty=G*u_\infty\,.
\end{equation}
Using the comparison principle in the radial variable, it has been proven in~\cite{bkln} that the radial, non-negative smooth solution to this problem is unique. In the original variables, the self-similar solutions of~\eqref{eq:KS} take the expression:
\begin{gather*}
\rho_\infty(t,x):=\frac 1{1+2t} u_\infty\left(\log(\sqrt{1+2t}), x/\sqrt{1+2t}\right)\,,\\
c_\infty(t,x):=v_\infty\left(\log(\sqrt{1+2t}),
x/\sqrt{1+2t}\right)\,.
\end{gather*}
%-----------------------------
\begin{theorem}[Large time behaviour,~\cite{BDP}]
 Under the assumptions in Theorem \ref{thm:Existence},
\[
\lim_{t\to\infty}\|{\rho(\cdot,t)-\rho_\infty(\cdot,t)}\|_{L^1(\R^2)}=0 \quad\mbox{and}\quad
\lim_{t\to\infty}\|{\nabla c(\cdot,t)-\nabla
c_\infty(\cdot,t)}\|_{L^2(\R^2)}=0\,.
\]
\end{theorem}
%-----------------------------
The proof follows the usual entropy/entropy production method in self-similar variables: We define the rescaled functions $u$ and
$v$ by:
\begin{equation*}
\rho(t,x)=\frac 1{R^2(t)} u\left(\frac x{R(t)},\tau(t)\right)
\quad\mbox{and}\quad c(t,x)=v\left(\frac x{R(t)},\tau(t)\right)
\end{equation*}
with $R(t)=\sqrt{1+2t}$ and $\tau(t)=\log R(t)$. The rescaled system is
\begin{equation}\label{Eqn:Keller-Segel-Rescaled}
\left\{
\begin{array}{ll}
\displaystyle \frac{\partial u}{\partial t}=\Delta u-\nabla \cdot (u(x+\nabla v))&\quad  x\in\mathbb{R}^2\,,\;t>0\,,\vspace{.1cm}\\
 \displaystyle v=G*u&\quad  x\in\mathbb{R}^2\,,\;t>0\,, \vspace{.1cm}\\
u(\cdot,t=0)=\rho_0 &\quad x\in\mathbb{R}^2\,,
\end{array}
\right.
\end{equation}
and the associated free energy takes the form
\begin{equation}\label{eq:dissip}
\F^R[u]:=\int_{\mathbb{R}^2}u \log u\dd x-\frac1{2}\int_{\mathbb{R}^2}u v\dd x +\frac12\int_{\mathbb{R}^2}|x|^2  u\dd x\,.
\end{equation}
If $(u,v)$ is a smooth solution of the rescaled Patlak-Keller-Segel system~\eqref{Eqn:Keller-Segel-Rescaled} which decays sufficiently at infinity, then
\[
\frac \dd{\dd t}\F^R[u(t,\cdot)]=-\int_{\mathbb{R}^2}u \left|\nabla \left(\log u -v+\frac{|x|^2}{2} \right)\right|^2  \dd x\,.
\]
If we keep in mind the gradient flow interpretation (which is true also for the rescaled equation), we can imagine that the limit when $t$ goes to infinity of a solution to the rescaled Patlak-Keller-Segel system~\eqref{Eqn:Keller-Segel-Rescaled} cancel the free energy dissipation~\eqref{eq:dissip}. So that the limit solution satisfies the Gelfand equation~\eqref{eq:stasol}.\medskip

The question of the speed of convergence is not fully understood. In the case of small mass, a linearisation method reduces the problem to a perturbation of the Fokker-Planck equation and gives:
%------------------------------------------------------------------------------
\begin{theorem}[Rate of convergence,~\cite{BDEF}]
 Under the assumptions in Theorem \ref{thm:Existence}, there exists a positive constant $M^*$ such that, for any initial data $\rho_0\in L^2(u_\infty^{-1}\dd x)$ of mass $M<M^*$, the rescaled Patlak-Keller-Segel system \eqref{Eqn:Keller-Segel-Rescaled} has a unique solution $u\in C^0(\R^+,L^1(\R^2))\cap L^\infty((\tau,\infty)\times\R^2)$ for any $\tau>0\,$. Moreover, there are two positive constants, $C$ and $\delta\,$, such that
\[
\int_{\R^2}|u(t,x)-u_\infty(x)|^2\, \frac{\dd x}{u_\infty(x)} \le C\,e^{-\,\delta\, t}\quad\forall\;t>0\;.
\]
As a function of $M\,$, $\delta$ is such that $\lim_{M\to 0_+}\delta(M)=1\,$.
\end{theorem}
%------------------------------------------------------------------------------

\quid{The question of the speed of convergence for larger mass remains open, even if numerical evidence indicates that it should be the case, see~\cite{BCC}.}
%%%%%%%%%%%%%%%%%%%%%%%%%%%%%%%%%
%%%%%%%%%%%%%%%%%%%%%%%%%%%%%%%%%
\subsection{Critical case}
  In the case $M=8\pi$, the free energy $\F$ is the same as the functional which appears in the logarithmic Hardy-Littlewood-Sobolev inequality~\eqref{eq:HLS}. The remainder entropy which was controlled in~\eqref{bdp} is thus entirely ``eaten'' by the logarithmic Hardy-Littlewood-Sobolev inequality~\eqref{eq:HLS}. In~\cite{BCM}, we use a three-step procedure:
  \begin{description}
  \item[How would it blowup] the space is split into balls and annulus. Using the diffusion it is possible to prove that in a ball the mass is less than $8\pi$ and to control the influence of the interactions outside the ball. So that only when all the mass is concentrated in a point, we cannot extend the solution to a bigger interval. If the solution blows up, it blows up as a Dirac mass concentrated in the centre of mass.
  \item[When would it blowup] in this case $M=8\pi$, by the virial computation~\eqref{eq:2moment}, the 2-moment remains constant. A De la Vall\'ee-Poussin's type argument, shows that the concentration cannot occur in finite time. If the solution blows up, it blows up as a Dirac mass concentrated in the centre of mass at infinite time.
  \item[Does it blowup] Keeping in mind the gradient flow structure described in the previous section, we can imagine that the solutions tend to the minimisers of the logarithmic Hardy-Littlewood-Sobolev inequality~\eqref{eq:HLS}. But as the second-moment is constant thanks to the virial computation~\eqref{eq:2moment}, the solution converges to the only minimiser $\rhohls$ of the logarithmic Hardy-Littlewood-Sobolev inequality~\eqref{eq:HLS} which is of finite moment: the Dirac mass.
  \end{description}
As a consequence, we prove
%------------------------------------------------------------------------------
\begin{theorem}[Infinite Time Aggregation, \cite{BCM}]\label{thm:main}
If the 2-moment is bounded, there exists a global in time non-negative
free-energy solution of the Patlak-Keller-Segel system~\eqref{eq:KS} with initial data $\rho_0$. Moreover
if $\{t_p\}_{p \in \N}\to\infty$ as $p\to\infty$, then $t_p
\mapsto \rho(t_p,x)$ converges to a Dirac of mass $8\,\pi$
concentrated at the centre of mass of the initial data weakly-* in the sense of measure as $p\to\infty$.
\end{theorem}
%------------------------------------------------------------------------------
In the radial and bounded case, the blowup rate and refined asymptotics estimates are given in the following
%---------------------------------------------------
\begin{theorem}[Blowup profile,~\cite{souplet}]
  In the radial case, in the ball, consider a solution to the Patlak-Keller-Segel system~\eqref{eq:KS} with $\partial \rho/\partial \nu=\rho \partial c /\partial \nu$ and $c=0$ on the boundary. Then when $t$ goes to $\infty$,
  \begin{equation*}
    \rho(t,0)=8 e^{5/2+2\sqrt{2}t}\left(1+O\left( t^{-1/2}\log(4t)\right) \right)\;.
  \end{equation*}
\end{theorem}
%---------------------------------------------------
In~\cite{SC02}, C. Sire and P.-H. Chavanis predicted this result by a formal argument considering only the first order correction terms. In~\cite{souplet}, N. Kavallaris and P. Souplet study the Patlak-Keller-Segel system~\eqref{eq:KS} on a bounded domain with Dirichlet condition, they make successive appropriate change of variable and of function, to reduce to a degenerate parabolic problem, so that the precise results are difficult to translate back to the original Patlak-Keller-Segel system~\eqref{eq:KS}. In particular, they prove that the solution is the sum of a quasi-stationary profile and of a correction term which is significant only for $x$ bounded away from $0$.
\medskip

The extension of Theorem~\ref{thm:main} to the case when the second moment is not finite allows the solution to converge to the other minimisers of the logarithmic Hardy-Littlewood-Sobolev inequality~\eqref{eq:HLS}. For this purpose we need to introduce another free energy functional which still has to be fully understood. Let us first recall the Fokker-Planck version of the fast diffusion equation corresponding to the fast diffusion equation ${\displaystyle \frac{\partial u}{\partial t} = \Delta \sqrt{u}}$ by a self-similar change of variable:
\begin{equation}\label{pmfp1}
\begin{cases}
 {\displaystyle \frac{\partial u}{\partial t}(t,x) =
 \Delta \sqrt{u(t,x)}  +
 2\sqrt{\frac{\pi}{\lambda M}} \,{\rm div}(x\,u(t,x))}\qquad & t>0\,,\;x\in\R^2\;,\vspace{.1cm}\\
 u(0,x)=u_0(x)\ge 0\qquad &x\in\R^2\, .
 \end{cases}
\end{equation}
This equation can also be written in a form analogous
to~\eqref{wasgrad}: following~\cite{lm} for $\lambda > 0$, define the relative entropy of the fast diffusion equation with respect to the stationary solution
$\rhohls$  by
\begin{equation*}%\label{Hlamdef}
 \Fcfd[u]:= \int_{\R^2} \frac{\left|\sqrt{u(x)}-\sqrt{\rhohls(x)}\right|^2}{\sqrt{\rhohls(x)}}\dd x
\end{equation*}
Equation~\eqref{pmfp1} can be rewritten as
\begin{equation*}
\frac{\partial u}{\partial t}(t,x) = {\rm
div}\left(u(t,x)\nabla \frac{\delta  \Fcfd[u(t)]}{\delta
u}(x)\right)\ ,
\end{equation*}
with
\begin{equation*}
\frac{\delta  \Fcfd[u]}{\delta u}  = \frac{1}{ \sqrt{\rhohls}} - \frac{1}{\sqrt{u}}\ ,
\end{equation*}

But the connection can be seen through the minimisers of $\Fcfd$ which are the same as those of the logarithmic Hardy-Littlewood-Sobolev inequality~\eqref{eq:HLS}. The functional $\Fcfd$ is a weighted distance between the solution and $\rhohls$ and its unique minimizer is $\rhohls$. It is tempting to compute the dissipation of $\Fcfd$ along the flow of solutions to the Patlak-Keller-Segel system~\eqref{eq:KS}: Let $\rho$ be a sufficiently smooth solution of the Patlak-Keller-Segel system~\eqref{eq:KS}. Then we compute
\begin{equation}\label{key1}
\frac{\dd}{\dd t} \Fcfd[\rho(t)]  =
-\frac{1}{2}\int_{\R^2}\frac{|\nabla \rho|^2}{\rho^{3/2}}\dd x +
\int_{\R^2}\rho^{3/2}\dd x  +
4\sqrt{\frac{M\,\pi}{\lambda}}\left(1- \frac{M}{8\pi}\right)\ .
\end{equation}
In the critical case $M=8\pi$ the dissipation of the $\Fcfd$ free energy along the flow of the Patlak-Keller-Segel system~\eqref{eq:KS} is
\begin{equation*}
  {\mathcal D}[\rho] := \frac{1}{2}\int_{\R^2}\frac{|\nabla \rho|^2}{\rho^{3/2}}\dd x - \int_{\R^2}\rho^{3/2}\dd x \ .
\end{equation*}
We use the following Gagliardo-Nirenberg-Sobolev inequality due to J.~Dolbeault and M.~Del Pino, see~\cite{DD}: For all functions $f$ in $\R^2$ with a square integrable
distributional gradient $\nabla f$,
\begin{equation*}
\pi \int_{\R^2} |f|^6 \dd x \le \int_{\R^2} |\nabla f|^2 \dd x \int_{\R^2}
|f|^4 \dd x \ ,
\end{equation*}
and there is equality if and only if $f$ is a multiple of a translate of $\rhohls^{1/4}$ for some $\lambda>0$.

As a consequence, taking $f = \rho^{1/4}$ so that $\int_{\R^2} f^4(x)\dd x = 8\pi$, we obtain ${\mathcal D}[\rho] \geq 0$, and moreover, ${\mathcal D}[\rho] = 0$ if and only $\rho$ is a translate of $\rhohls$ for some $\lambda>0$.
%----------------------
\begin{remark}\label{benasque} While I was finishing this survey, J. A. Carrillo and I realised that this new free energy gives another proof of non existence of global-in-time solutions in the super-critical case. Indeed, by ${\mathcal D}[\rho] \ge 0$ and~\eqref{key1}, $0 \le \Fcfd[\rho(t)] \le 4\sqrt{{M\,\pi}/{\lambda}}\left(1- {M}/{8\pi}\right) t$. So that in the case $M>8\pi$, if $\Fcfd[\rho_0]$ is initially bounded, and thus of infinite 2-moment, then there cannot be global-in-time solutions.  
\end{remark}
%--------------------
\medskip

Based on this remark, the main results of \cite{BCCxx} in the critical case can be summarised in the following:
%--------------------------------------------------------
\begin{theorem}[Existence of global solutions,~\cite{BCCxx}]\label{main}
Let $\rho_0$ be any density in $\R^2$ with mass $8\pi$, such that  $\Fhls[\rho] < \infty$, and for some $\lambda>0$,
$\Fcfd[\rho] < \infty$.  
Then there exists a global free energy solution of the Patlak-Keller-Segel
equation~\eqref{eq:sp} with initial data $\rho_0$.
Moreover,
$$\lim_{t\to\infty}\Fhls[\rho(t)] = \Fhls[\rhohls]\qquad{\rm and}\qquad   \lim_{t\to\infty}\|\rho(t) - \rhohls\|_1 = 0\ .
$$
\end{theorem}
%--------------------------------------------------------
We even prove further regularity using the propagation of the $L^p$-estimates and the hyper\-contractivity property of the equation. This theorem can be translated in terms of the Wasserstein distance of the initial data to $\rhohls$ thanks to the Talagrand inequality:
\begin{equation*}
\W_2(\rho,\rhohls) \le\sqrt{\frac{2\Fcfd[\rho]}{2\sqrt{\frac{\pi}{M\,\lambda}}}}.
\end{equation*}
If we keep in mind that the 2-moment can be seen as the Wasserstein distance between the solution and the Dirac mass, we see that Theorem~\ref{main} completes the picture which emerged from the first result set out at the beginning of the section. As soon as we start at a finite distance from one of the minimisers $\rhohls$ we can construct a solution which converges towards it. Note that this result is true for the solution that we construct as we do not have uniqueness of the solution to the Patlak-Keller-Segel system, even if we strongly believe that this is the case. Also observe that the equilibrium solutions $\rhohls$ are infinitely far apart: let
$\varphi(x) = \sqrt{\lambda/\mu}|x|^2/2$, one has $\nabla\varphi\#\varrho_{\mu} = \rhohls$. Thus, 
$$
\W_2^2(\varrho_{\mu} ,\rhohls ) =
\frac{1}{2}\int_{\R^2}\left| \sqrt{\frac{\lambda}{\mu}} x
-x\right|^2 \varrho_{\mu}(x) \dd x = +\infty
$$
since the equilibrium densities $\rhohls$ all have
infinite second moments. In particular,
$\Fcfd[\varrho_{\mu}]=+\infty$ for $\mu\neq\lambda$. There may still be initial data out of these basins of attraction.\medskip

Concerning the proof of Theorem~\ref{main}, we expect the propagation of the bounds on $\Fhls[\rho]$ and ${\mathcal D}[\rho]$ to give compactness. Unfortunately, ${\mathcal D}[\rho]$ is a difference of two functionals of $\rho$ that can each be arbitrarily large even when
${\mathcal D}[\rho]$ is very close to zero. Indeed, for $M = 8\pi$
and each $\lambda>0$, ${\mathcal D}[\rhohls] = 0$ while
$$
\lim_{\lambda\to 0}\|\rhohls\|_{3/2} = \infty \, , \quad
\lim_{\lambda\to 0}\|\nabla\rhohls^{1/4}\|_{2} = \infty \, , \quad
\mbox{and} \quad \lim_{\lambda\to 0}\rhohls = 8\pi\delta_0\, .
$$
Likewise, an upper bound on $\Fhls[\rho]$  provides no upper bound on the entropy $\int_{\R^2}\rho\log\rho$. Indeed,
$\Fhls[\rho]$ takes its minimum value for $\rho= \rhohls$ for each
$\lambda > 0$, while  $\lim_{\lambda\to 0}\int \rhohls \log \rhohls =
\infty$. Fortunately, an upper bound on both $\Fcfd[\rho]$ and
$\Fhls[\rho]$ {\em does} provide an upper bound on $\int \rho \log \rho$:
%------------------------------------
\begin{theorem}[Concentration control for $\Fhls$, \cite{BCCxx}]\label{both}
Let $\rho$ be any density with mass $M = 8\pi$ such that $\Fcfd[\rho] < \infty$
for some $\lambda > 0$.  Then there exist $\gamma_1>0$
and $C_{{\rm CCF}}>0$ depending only on $\lambda$ and $\Fcfd[\rho]$
such that 
\begin{equation*}
\gamma_1\int_{\R^2}\rho \log \rho \dd x  \le \Fhls[\rho]  +
C_{{\rm CCF}}\,.
\end{equation*}
\end{theorem}
%------------------------------------
Here we also prove that since $\Fcfd$ controls concentration, a uniform bound on both $\Fcfd$ and ${\mathcal D}$ does indeed provide compactness. We shall prove:
%------------------------------------
\begin{theorem}[Concentration control for ${\mathcal D}$, \cite{BCCxx}]\label{both2}
Let $\rho$ be any density in $L^{3/2}(\R^2)$ with mass $8\pi$ such that $\Fhls[\rho]$ is finite, and $\Fcfd[\rho]$ is finite for some $\lambda> 0$. Then there exist constants $\gamma_1>0$ and $C_{{\rm CCF}}>0$ depending only on $\lambda$, $\Fcfd[\rho]$ and $\Fhls[\rho]$ such that
\begin{equation*}
\gamma_2\,\int_{\R^2}|\nabla \rho^{1/4}|^2 \dd x \le \pi{\mathcal
D}[\rho] +C_{{\rm CCD}}\ \, .
\end{equation*} 
\end{theorem}
%------------------------------------
The proofs of this two theorems leads on the following lemma:
\begin{equation}\label{1mom}
  \int_{\R^2}\sqrt{\lambda +|x|^2}\,\rho(x)\dd x \le 2\,\sqrt\lambda\,M + 2M^{3/4}(\lambda/\pi)^{1/4}\,\sqrt{\Fcfd[\rho]}\;.
\end{equation}
%------------------------------------

As explained at the beginning of the section, in~\cite{BCM} we managed to find a ball in which the mass was smaller than $8\pi$. Here,~\eqref{1mom} gives a vertical cut to prove Theorem~\ref{both}. Indeed, we split the function $\rho$ in two parts: given $\beta > 0$, define $\rho_{\beta}(x)  = \min\{ \rho(x)\ , \ \beta\}$. By~\eqref{1mom}, for $\beta$ large enough, $\rho-\rho_{\beta}$ is such that:
\begin{equation*}
\int_{\RR^2}\left(\rho-\rho_\beta\right) \le
\frac{C_1}{\beta} + C_2\,\sqrt{\Fcfd[\rho]}  \le
\frac{C_1}{\beta} + \frac{8\pi - \eps_0}{2} < 8\,\pi -
\varepsilon_0\;.
\end{equation*}
We then apply the logarithmic Hardy-Littlewood-Sobolev inequality method as in~\eqref{eq:entro} to the function $\rho-\rho_{\beta}$ whose mass is less than $8\pi$.

The same idea works for the Gagliardo-Nirenberg-Sobolev inequality to prove Theorem~\ref{both2}: Let $f := \rho^{1/4}$, we split $f$ in two parts by defining $f_\beta:=\min\{f,\beta^{1/4}\}$ and $h_\beta:=f-f_\beta$. We use~\eqref{1mom} and apply the Gagliardo-Nirenberg-Sobolev inequality to control $h_\beta$.\medskip

\noindent{\bf Idea of the proof of Theorem~\ref{main}:} It follows the line of the convergence of the JKO minimising scheme~\eqref{JKOscheme} exposed in the previous section to obtain the Euler-Lagrange equation~\eqref{eq:EL}. Dividing the Euler-Lagrange equation~\eqref{eq:EL} by $\sqrt {\rho_\tau^{n+1}}$ we obtain:
\begin{equation}\label{eq:eldicrho}
2\nabla \sqrt{\rho_\tau^{n+1}} = \left(\nabla c_\tau^{n+1} - \frac{x-\nabla \varphi_\tau^{n}}{\tau}\right)\sqrt{\rho_\tau^{n+1}}\;,
\end{equation} 
where $\nabla \varphi_\tau^{n} \# \rho_\tau^{n+1} = \rho_\tau^{n}$. Integrating~\eqref{eq:eldicrho} we obtain
\begin{equation*}
    \int_{\R^2} \left|\sqrt{\rho_\tau^{n+1}}\nabla c_\tau^{n+1} - 2\nabla \sqrt{\rho_\tau^{n+1}}\right|^2\dd x=\int_{\R^2}\left|\frac{x-\nabla \varphi_\tau^{n}}{\tau}\right|^2\rho_\tau^{n+1}\dd x = \W_2(\rho_\tau^{n},\rho_\tau^{n+1})
\end{equation*}
which is bounded thanks to~\eqref{eq:estime}. But the left hand side is a sum of two terms so that we cannot conclude any compactness on each of them. That is the reason why we have to introduce a regularisation of the Green kernel:
let $\gamma$ be the standard Gaussian probability density in $\R^2$. Then, for all $\epsilon>0$ define $\gamma_\epsilon(x) =\epsilon^{-2}\gamma\left(x/\epsilon\right)$, and define  the monotonically regularised Green function $G_\epsilon = \gamma_\epsilon*G*\gamma_\epsilon$. And the corresponding regularised free energy functional defined, for all $0<\epsilon\le1$, by
\begin{equation*}
\Fhls^{\epsilon}[\rho] := \int_{\R^2}\rho(x) \log \rho (x)\dd x -
\frac{1}{2}\iint_{\R^2 \times\R^2} \rho(x)\,G_\epsilon(x-y)\,\rho(y)\,\dd x \dd y\;.
\end{equation*}
The regularised chemical attractant density $c_\epsilon:=G_\epsilon *\rho(x)$ is uniformly bounded, so that by the triangle inequality, from~\eqref{eq:eldicrho} we obtain
\begin{align*}
2\|\nabla \sqrt{\rho}\|_2 &\le  \left(\int_{\R^2}\left|\nabla c_\epsilon(x)\right|^2\rho(x)\dd x\right)^{1/2}
+\frac1{\tau} \left(\int_{\R^2}|x-\nabla \varphi(x)|^2\rho(x)\dd x\right)^{1/2}\nonumber\\
&\le  \left(\int_{\R^2}\left|\nabla
c_\epsilon(x)\right|^2\rho(x)\dd x\right)^{1/2} + \frac1{\tau}
\W_2(\rho,\rho_0)\ .
\end{align*}
Which gives a bound, which depends on $\eps$, on $\|\nabla \sqrt{\rho}\|_2$ and so on $\|\rho\|_p$ for all $p \in (1,\infty)$ by the Gagliardo-Nirenberg-Sobolev inequality~\eqref{GNS}. The uniform estimates are of course more difficult to get. We could probably prove the existence by using the regularisation procedure but it is more elegant to prove it following the JKO minimisation scheme~\eqref{JKOscheme}. This optimal transport framework is better adapted because in this case we can use displacement convexity. The displacement convexity of $\Fcfd$ is formally obvious from the fact that
\begin{equation*}
 \Fcfd[u]=\int_{\R^2} \left( -2 \sqrt{u(x)} +\sqrt{\frac{1}{2\lambda}} \frac{|x|^2}{2}u(x)  \right)\dd x+C \;.
\end{equation*}
where $-\sqrt{u(x)}$ and $|x|^2u(x)$ are displacement convex.  Actually, we can not expand the terms as none of them are finite. Thus, we are forced to introduce a regularisation of $\Fcfd$.  While there are many tools available to regularise functions that are convex in the usual sense, there is no general approach to regularising functionals while preserving, or at least not severely damaging, their formal displacement convexity properties. The following regularisation is one of the cornerstones of this proof:
\begin{equation*}
\Fcfdr[u] = \int_{\R^2}\frac {\left(\sqrt{u+\delta } - \sqrt{
\rhohls+\delta }\right)^2} {\sqrt
{\rhohls+\delta} }\dd x\
\end{equation*}
which is $\gamma_\delta$-displacement convex with $\gamma_\delta \le 0$.

At each step, the main {\it a priori} estimate comes from the convexity estimate of the type
\begin{equation*}
  \Fcfd[\rho_0] -\Fcfd[\rho_1] \ge \limsup_{t \to 0}\frac{\Fcfd[\rho_t]-\Fcfd[\rho_0]}{t}\;,
\end{equation*}
which give in this optimal transport framework the {\em above the tangent formulation}: 
\begin{equation*}
 \Fcfd[\rho_\tau^{n+1}] - \Fcfd[\rho_\tau^{n}]   \ge     \frac12\int_{\R^2} \left[\sqrt{\frac{1}{2\lambda}}\,x +
 \frac{\nabla \rho_\tau^{n}}{\left(\rho_\tau^{n}\right)^{3/2}}\right] \cdot (\nabla \psi(x) -x)\,\rho_\tau^{n}\dd x\ .
\end{equation*}
where $\nabla \psi$ is such that $\nabla \psi\# \rho_\tau^{n+1} = \rho_\tau^{n}$. Inferring the Euler-Lagrange equation~\eqref{eq:EL}: $-\nabla \rho_\tau^{n+1} + \rho_\tau^{n+1} \nabla c_\tau^{n+1} = \left({\rm id} - \nabla \psi\right)\rho_\tau^{n+1}/\tau$, we obtain a discrete version of the entropy/entropy dissipation inequality
\begin{equation}\label{eq:convexest}
\Fcfd[\rho_\tau^{n+1}] - \Fcfd[\rho_\tau^{n}]\le -\tau\,{\mathcal D}[\rho_\tau^{n}]\ .
\end{equation}
This inequality is a skeleton version of the crucial estimate which allow to apply the standard entropy/entropy dissipation method to study the asymptotics.\medskip

We have to pass to the limit in three parameters: $\delta$, $\eps$ and $\tau$. The $\delta$-regularisation was only introduced to justify the integration by parts in~\eqref{eq:convexest}. While passing to the limit in $\delta$ we lose slightly in the convexity estimate~\eqref{eq:convexest}. We obtain that there exist $\AA$ such that
\begin{equation}\label{eq:witouteps}
\Fcfd[\rho] \le \Fcfd[\rho_0]-\tau\,{\mathcal D}[\rho] +\tau\,\AA\,\|\gamma\|_{4/3}
\quad\mbox{and}\quad
\Fcfd[\rho] \leq \Fcfd[\rho_0]-\tau \,{\mathcal D}[\rho] +
\tau\sqrt\epsilon\,\AA\,\| \rho \|_{4/3}\ .
\end{equation}
However, a tricky interplay between $\eps$ and $\tau$ allows to conserve the estimate at the limit. Let us give an idea on how to do it at the first step: define $Q_0 > 0$, $\tau_0^{\star}>0$ by
\begin{equation}\label{taude}
Q_0:=C_{\rho_0}-  \Fcfd[\rho_0]\quad \mbox{and} \quad
  \tau_0^\star := \min\left\{\frac{Q_0}{2\,\AA\|\gamma\|_{4/3}}\ ,\ 1\right\}\ ,
\end{equation}
and then $\epsilon_\ell$ which depends on $\tau$ by
\begin{equation}\label{taude2}
 \tau^{1/3}\sqrt{\epsilon_\ell}\left[8\,\pi^{1/3}\,\AA\, \gamma_2^{-2/3} \left(\pi \,C_{\rho_0}+\tau_0^\star C_{{\rm CCD}} \right)^{2/3}\right] = \frac{Q_0}{4}\tau^2 2^{-\ell}\ .
\end{equation}
By~\eqref{eq:witouteps}, our choice of $\tau$ and $Q_0$ in~\eqref{taude} implies that
\begin{equation}\label{fonda}
  \Fcfd[\rho] \le \Fcfd[\rho_0] -\tau{\mathcal D}[\rho]+ \frac{Q_0}{2} =
  C_{\rho_0} - Q_0 -\tau{\mathcal D}[\rho]+ \frac{Q_0}{2}\le C_{\rho_0}-\tau{\mathcal D}[\rho] \ .
\end{equation}
On one hand, the Gagliardo-Nirenberg-Sobolev inequality~\eqref{GNS} induces ${\mathcal D}[\rho] \ge 0$ so that~\eqref{fonda} implies that the energy estimates on $\rho$ needed for the two concentration controlled inequalities, Theorems~\ref{both} and~\ref{both2} propagate in the sense that
\begin{equation*}
\Fhls[\rho] < +\infty\;, \quad \Fcfd[\rho]  < C_{\rho_0} \;.%adrien
\end{equation*}
On the other hand, since $\Fcfd[\rho]$ cannot be negative it implies
\begin{equation*}
{\mathcal D}[\rho] \le \frac{C_{\rho_0}}{\tau}\ .
\end{equation*}
The concentration controlled inequality, Theorem~\ref{both2}, thus gives a bound on $\|\nabla \rho^{1/4}\|_2$:
$$
 \int_{\R^2} \left|\nabla \rho^{1/4}\right|^2 \dd x \le
 \frac{1}{\gamma_2}\left[\pi{\mathcal D}[\rho] + C_{{\rm CCD}}\right] \le
 \frac{1}{\tau} \frac{1}{\gamma_2}\left[\pi C_{\rho_0} + \tau_0^\star C_{{\rm CCD}}\right] \ .
$$
By the Gagliardo-Nirenberg-Sobolev inequality we obtain a bound of the form
\begin{equation*}
\int_{\R^2}\rho^{3/2} \dd x \le 8 \int_{\R^2} \left|\nabla
\rho^{1/4}\right|^2 \dd x \le  \frac{1}{\tau}
\frac8{\gamma_2}\left[\pi C_{\rho_0}+\tau_0^\star C_{{\rm CCD}}
\right] := \frac{C_3}{\tau} \ .
\end{equation*}
And, by H\"older's inequality,
\begin{equation*}%\label{hold}
\int_{\R^2}\rho^{4/3} \dd x = \int_{\R^2}\rho^{1/3} \rho \dd x \le
(8\pi)^{1/3}\left(\int_{\R^2}\rho^{3/2} \dd x\right)^{2/3}\le
(8\pi)^{1/3}\left( \frac{C_3}{\tau} \right)^{2/3}\ .
\end{equation*}
Now using this bound in~\eqref{eq:witouteps}, we obtain
\begin{equation*}
\Fcfd[\rho] -\Fcfd[\rho_0] \le -\tau{\mathcal D}[\rho]
+\tau^{1/3}\sqrt{\epsilon}\,\left[\AA\, (8\pi)^{1/3} C_3^{2/3}\right]
-2\sqrt{\frac{\pi}{M\,\lambda}} \W_2^2(\rho_0,\rho).
\end{equation*}
Our choices of $\eps$ in~\eqref{taude2} ensures that we do not lose too much in the convexity estimate:
\begin{equation*}
\Fcfd[\rho] -\Fcfd[\rho_0] \le -\,\tau{\mathcal D}[\rho] +
\frac{Q_0}{4} \tau^22^{-\ell} \ .
\end{equation*}
We see here that the extra term ${Q_0} \tau^22^{-\ell}/{4}$ in the convexity estimate remains bounded in the sum. It actually converges to $0$ when $\tau$ goes to $0$. These proofs are technical and, even if some are interesting from the general theoretical point of view and could be used for other problems. Moreover, we can propagate the $L^p$-estimates and prove hyper-contractivity and to obtain the convergence of the JKO minimising scheme~\eqref{JKOscheme} toward a free energy solution to the Patlak-Keller-Segel system~\eqref{eq:KS}. The proof of the asymptotics result is then a standard entropy method as set out at the end of Section~\ref{subcritical}.

Very recently, in~\cite{CFxx}, E. Carlen and A. Figalli use a argument of Bianchi-Egnell to obtain a quantitative stability for the logarithmic Hardy-Littlewood-Sobolev inequality~\eqref{eq:HLS} and prove:
\begin{equation*}
  \|\rho(t)-\rhohls \|_{L^1(\R^2)} \le \frac{C}{\sqrt{\log(e+t)}}\;.
\end{equation*}
%%%%%%%%%%%%%%%%%%%%%%%%%%%%%%%%
\section{The Non-linear parabolic-elliptic Patlak-Keller-Segel system}
\subsection{The model}
In higher dimensions the critical quantity is no longer the mass but the $L^{d/2}$-norm, see~\cite{CPZ}. We can however replace the linear diffusion with a homogeneous non-linear diffusion:
\begin{equation}\label{eq:sp}
\left\lbrace
\begin{array}{rll}
\displaystyle \frac{\partial \rho}{\partial t}(t,x)&={\rm div } \left[
\nabla \rho^m(t,x)-\rho(t,x)\nabla \phi(t,x)\right]\qquad &
t>0\,,\;x\in\RR^d\;,\vspace{.3cm}\\
\displaystyle -\Delta \phi(t,x)&=\rho(t,x)\;,\qquad
&t>0\,,\;x\in\RR^d\;,\vspace{.3cm}\\
\rho(0,x)&=\rho_0(x)\qquad &x\in\RR^d\, ,
\end{array}
\right.
\end{equation}
where $m \in (0,1)$ and $d \ge 3$. In astro-physics, this system models the motion of the mean field of many self-gravitating Brownian particles, including a version of the Chandrasekhar equation modelling gravitational equilibrium of poly-tropic stars~\cite{Chandrasekhar}. This system is then known as the generalised Smulochowski-Poisson system, in dimension $3$ see~\cite{MR1052381} and~\cite{biler-nadzieja, ChS04} in dimension $1$ and $2$.

Define
\begin{equation*}
\EE(x)=c_d \;\frac{1}{|x|^{d-2}} \quad\mbox{and}\quad
  c_d:=\frac {1}{(d-2) \sigma_d}
\end{equation*}
where $\sigma_d:=2\,\pi^{d/2}/\Gamma(d/2)$ is the surface area of the
sphere $\mathbb{S}^{d-1}$ in $\RR^d$. Up to a harmonic function $\phi=\EE\ast \rho$, so that the system~\eqref{eq:sp} can be rewritten as a non-local parabolic equation:
\begin{equation}\label{eq:spf}
\frac{\partial \rho}{\partial t}(t,x)={\rm div }
\left[ \nabla \rho^{m}(t,x)-\rho(t,x)
\nabla (\EE \ast \rho)(t,x)\right]\qquad  t>0\,,\;x\in\RR^d\, .
\end{equation}
In this case too the mass is preserved and will be denoted $M$.

Let $\rho_\lambda(x):=\lambda^d \rho(\lambda\,x)$ with $\lambda>0$, the diffusion term scales like $\lambda^{d\,m+2}\Delta (\rho_\lambda^m)(\lambda\,x)$ whereas the interaction term scales like $\lambda^{2\,d}{\rm div }\left(\rho_\lambda\nabla (\EE \ast \rho_\lambda)\right)(\lambda\,x) $. Hence the mass-invariant scaling of the diffusion term balances the potential drift in \eqref{eq:spf} if
\begin{equation}
  \label{eq:md}
  m=m_d=:2 \left(1 - \frac{1}{d}\right)\in (1,2)\,.
\end{equation}
This difference of balance was studied to obtain
%--------------------
\begin{theorem}[First criticality,~\cite{S1,S2}]
Let $m_d$ be as defined in~\eqref{eq:md}.
  \begin{itemize}
  \item if {$m>m_d$} then the solutions to~\eqref{eq:sp} exist globally in time,
  \item if {$m<m_d$} then solutions to~\eqref{eq:sp} with sufficiently large initial data blowup in finite time,
  \item if {$m=m_d$} there exist two constants {$M_1>0$} and {$M_2>M_1$} such that
    \begin{itemize}
    \item if {$M<M_1$} then the solutions to~\eqref{eq:sp} exist globally in time,
    \item if {$M>M_2$} there exist initial conditions such that the corresponding solution blows up in finite time.
    \end{itemize}
  \end{itemize}
\end{theorem}
%--------------------
The proof of this theorem relies on the Gagliardo-Nirenberg-Sobolev inequality and is not sharp as we will see in Theorem~\ref{sharp} below. In dimension 2, corresponding results were obtained by~\cite{Kowalczyk04,calvezcarrillo}. When no confusion is possible, the index $d$ in $m_d$ will be omitted and the critical exponent will be denoted $m$ in the sequel of this article. 
\medskip

The analogous of the free energy used in the previous section is:
\begin{equation*}
 t\mapsto \GG[\rho(t)]:= \int_{\RR^d} \frac{\rho^m(t,x)}{m-1} -
 \frac12 \iint_{\RR^d \times \RR^d} \EE(x-y)\,\rho(t,x)\,\rho(t,y)\dd x \dd y\;
\end{equation*}
which is related to its time derivative along the flow of~\eqref{eq:spf} by
\begin{equation*}
\frac\dd{\dd t}\GG[\rho(t)]=-\int_{\RR^d}\rho(t,x)\left|\nabla\left(\frac{m}{m-1}
\rho^{m-1}(t,x) -\phi(t,x) \right)\right|^2\dd x \; .
\end{equation*}
%%%%%%%%%%%%%%%%%%%%%%%%%
\subsection{The sub-critical case}
In~\cite{BCL09}, the functional inequality used is a {\it variant to the
Hardy-Littlewood-Sobolev (VHLS) inequality}: for all $h\in \LL^1(\RR^d)\cap \LL^m(\R^d)$,
there exists an optimal constant $C_*$ such that
\begin{equation*}
C_*=\sup_{h\neq0}\left\{ \|h\|^{-m}_{m}\,\|h\|^{-{2}/{d}}_{1} \iint_{\RR^d\times\RR^d}\frac {h(x)\,h(y)}{|x-y|^{d-2}} \dd
x\dd y\right\}\;.
\end{equation*}
We define the critical mass by 
\begin{equation*}
  M_c:=\left[\frac{2}{(m-1)C_*c_d} \right]\;.
\end{equation*}
%---------------------------------
\begin{theorem}[Global-in-time existence,~\cite{BCL09,ST09}]\label{sharp}
  If $u_0$ is of mass $M<M_c$ then there exists a global weak solution with initial condition $u_0$. Moreover, this solution satisfies the free energy/free energy dissipation inequality.
\end{theorem}
%---------------------------------
The proof of existence follows the lines of the 2d Patlak-Keller-Segel system. Indeed, as a direct consequence of the VHLS inequality, for any solution $\rho$ to the nonlinear Patlak-Keller-Segel system~\eqref{eq:spf}
\begin{equation*}
\frac{C_*\,c_d}{2}\, \left(M_c^{2/d}  - M^{2/d}\right)\|\rho(t)\|_m^m \le \GG[\rho(t)]\le \GG[\rho_0] <\infty\,.
\end{equation*}
In the case $M<M_c$, it gives the concentration controlled analogous to the entropy {\it a priori} estimate~\eqref{eq:entro} of the previous section. 
It should not be difficult to prove the existence of global-in-time solutions using the JKO minimising scheme.\medskip

\quid{The convergence to the self-similar solution is still open. By doing the porous medium scaling, we can prove, see~\cite[Theorem~5.2]{BCL09}, that for any given mass $M<M_c$ there exists a unique minimiser $W_M$ to the rescaled free energy. Moreover this minimiser is non-negative, radially symmetric and compactly supported. We expect this minimiser to attract all the solutions but we have not been able to prove it.}
%%%%%%%%%%%%%%%%%%%%%%%%%
\subsection{The critical case}
The balance in the mass-invariant scaling of diffusion and potential drift can also be seen in the free energy: If $h_\lambda(x):=\lambda^dh(\lambda\,x)$ then
\begin{equation*}
  \mathcal G[h_\lambda]=\lambda^{{(m-1)d}}\int_{\R^d} \frac{h^m(x)}{m-1} - \lambda^{{d-2}}\frac{c_d}{2} \iint_{\R^d \times \R^d}\frac{1}{|x-y|^{d-2}}\,h(x)\,h(y)\dd x \dd y
\end{equation*}
The diffusion and interaction term balance if $m=m_d$. And in this case 
\begin{equation}\label{scaleg}
  \mathcal G[h_\lambda]={\lambda^{d-2}\mathcal G[h]}\;.
\end{equation}
This scaling of $\mathcal G$ has to be compared to the case $d=2$, see~\eqref{eq:scalingf}. This major difference will be a key in the application of the concentration-compactness argument. In the case of the 2d Patlak-Keller-Segel model~\eqref{eq:KS}, the authors were not able to apply a concentration compactness argument due to the rigidity in the scaling of the free energy. Indeed, the scaling which preserves the mass also preserves the free energy. Here, as will be shown below, it is possible to follow the line of Pierre-Louis Lions's original article.\medskip

The minimisers of $\GG$ of mass $M_c$ are such that there are $R>0$ and $z \in \RR^d$ with
  \begin{equation}\label{eq:defv}
   V(x)= \left\{
      \begin{array}{ll}
        \displaystyle \frac{1}{R^d} \left[\zeta\left(\frac{x-z}{R}
\right) \right]^{d/(d-2)} \quad & \mbox{if $x \in B(z,R)$,}\vspace{.3cm}\\
\displaystyle 0 \quad & \mbox{if $x \in \RR^d\setminus B(z,R)$}
      \end{array}
 \right.
  \end{equation}
where $\zeta$ is the unique positive radial classical solution to
$$
\Delta \zeta + \frac{m-1}{m}\, \zeta^{1/(m-1)}=0 \;\;\mbox{ in }\;\;
B(0,1) \;\;\mbox{ with }\;\; \zeta=0 \;\;\mbox{ on }\;\;\partial
B(0,1)\,.
$$

\quid{On the contrary to the 2d Patlak-Keller-Segel system, this proves that there exist compactly supported stationary solutions not blowing-up at infinite time. We were however not able to prove that they attract some solutions.}

Here, thanks to~\eqref{scaleg} we can adapt the concentration-compactness method to prove
%-------------------------------------------------
  \begin{proposition}[How would it blowup,~\cite{BCL09}]
    Let $T \in (0,\infty]$ and a sequence $(t_k)_k$ converging to $T$. If 
{$$
\lim_{k \to \infty}\|\rho(t_k)\|_m=\infty\;.
$$}
then there are a sub-sequence $(t_{k_j})_j$ and a sequence $(x_j)_j$ in {$\RR^d$} such that
\begin{equation*}
  \lim_{j \to \infty} \left\| \rho(t_{k_j},x+x_j)-\frac{1}{\lambda_{k_j}^d}\, {V}\left(\frac{x}{\lambda_{k_j}} \right)\right\|_{L^1}=0\;,
\end{equation*}
where {$\lambda_k:=\|\rho(t_k)\|_m^{-m/(d-2)}$} and {$V$ is the minimiser of $\mathcal G$ of the form~\eqref{eq:defv} with $\|V\|_m=1$}.
  \end{proposition}
%-------------------------------------------------
The main ingredient of the proof is the following: we set $v_k(x):=\lambda_k^d\, \rho(t_k,\lambda_k\,x)$ so that {$\|v_k\|_m=1$}. By the concentration compactness principle there exists a sub-sequence satisfying compactness, vanishing or dichotomy. As already discussed, contrary to the 2d Keller-Segel system, here 
\begin{equation*}
 \lim_{k \to \infty} \mathcal G[v_k]=\lim_{k \to\infty}{\|\rho(t_k)\|_m^{-m}}\,\GG[\rho(t_k)] \;{ = 0}\;. 
\end{equation*}
As a consequence
\begin{equation*}
   \lim_{k \to \infty}\iint_{\R^d \times \R^d}\frac{v_k(t,x)\,v_k(t,y)}{|x-y|^{d-2}}=  \lim_{k \to \infty}
\frac{2}{c_d}\left(\frac{1}{m-1} \|v_k\|_m^m -  \GG[v_k] \right)>0 \;. 
\end{equation*}
Whereas this quantity goes to zero if vanishing or dichotomy should occur.
%
%\quid{This theorem gives the profile of the blowup if the solutions blowsup but we do not know if such critical mass blowing up solution exists.}

\quid{Except in the radially symmetric case, we are not able to say if the blowup occurs at the centre of mass or if the blowup escapes at infinity. We could not even rule out the possibility that a sub-sequence diverges whereas the other does not.}
%%%%%%%%%%%%%%%%%%%%%%%%%
\subsection{The super-critical case}
The answer was clear for the super-critical case for the 2d Patlak-Keller-Segel thanks to the constant sign of the derivative of the 2-moment. Here for any $\rho$ solution to~\eqref{eq:sp} the virial identity is
 \begin{equation}\label{virialm}
  \frac{\dd}{\dd t} \int_{\R^d}|x|^2\,\rho(t,x)\dd x = 2\,(d-2)\,{\mathcal G[\rho(t)]}\;.
 \end{equation}
 \begin{theorem}[Blowup,~\cite{S1,S2,BCL09,ST09}]
   If $M>M_c$, there exist initial data of mass $M$ such that the $L^m$-norm of the corresponding solution blows up in finite time. 
 \end{theorem}
In~\cite{BCL09}, the proof relies on a procedure directly adapted from~\cite{Weinstein83}: let $\tilde{\rho}$ be a minimiser of the form~\eqref{eq:defv} and consider
\begin{equation*}
  \rho_0=\frac{M}{M_c}\tilde{\rho}\;.
\end{equation*}
Then, 
\begin{equation*}
  \mathcal G[\rho_0]=\frac{1}{m-1}\left(\frac{M}{M_c} \right)^m \left[ 1-
\left({\frac{M}{M_c}} \right)^{2-m}\right]\|\tilde{\rho}\|_m^m
\end{equation*}
is negative if $M \ge M_c$. This result combined with~\eqref{virialm} gives the expected result.

\quid{We cannot exclude the possibility that solutions with positive free energy exist globally in time.}
\medskip

\begin{figure}[h!]
 \begin{minipage}[t]{1\linewidth}
\centering\epsfig{figure=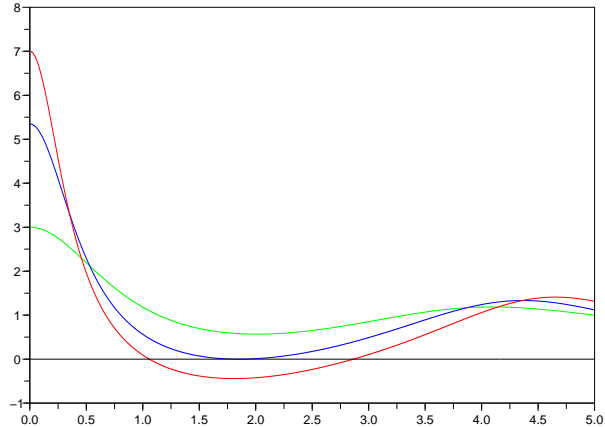,width=7cm,angle=270}
  \caption{Behaviour of $u(.,a)$ for $a>a_c$, $a=a_c$ and $a<a_c$.}\label{fig:lemma5}
 \end{minipage} \hfill
\end{figure}

In~\cite{bl} a more precise partial answer is proven: the blowup time $T$ being given, we can look for solution to~\eqref{eq:sp} of the form
\begin{equation*}
  \rho(t,x)=\frac{1}{s(t)^d}\Psi \left( \frac{x}{s(t)}\right)\quad \mbox{and}\quad c(t,x)=\frac{1}{s(t)^{d-2}}\Phi \left( \frac{x}{s(t)}\right)
\end{equation*}
where $s(t):=[d(T-t)]^{1/d}$.
%------------------------------
\begin{theorem}[Self-similar blowing-up solutions,~\cite{bl}]\label{selfsi}
There exists {$\Tilde M_c \in (M_c,\infty)$} such that for any $M \in (M_c,\Tilde M_c]$, there exists a self-similar blowing-up solution with a radially symmetric, compactly supported {and non-increasing} profile $\Psi$, satisfying $\|\rho(t)\|_1=\|\Psi\|_1=M$ for $t \in [0,T)$ and $\|\rho(t)\|_\infty \to \infty$ as $t \to T$.  
\end{theorem}
%------------------------------
The method relies on the study of a boundary value problem for the following non-linear ordinary differential equation:
\begin{equation*}
\left\{
\begin{array}{l}
\displaystyle{u''(r,a) + \frac{d-1}{r}\ u'(r,a) + |u(r,a)|^{p-1}\ u(r,a) - 1 = 0}\,, \quad r\in [0,r_{\max}(a))\,,\vspace{.1cm}\\
u(0,a)= a\,, \quad u'(0,a)=0\,,
\end{array}
\right.
\end{equation*}
with $r_{\max}(a)\in (0,\infty]$ and $p=d/(d-2)$. We prove that there are global solutions to this problem and that the solutions oscillate around the stationary solution 1, see Figure~\ref{fig:lemma5}. The solution of Theorem~\ref{selfsi} corresponds to that which vanishes and the support corresponds to its first zero. The proofs rely on ordinary differential equation tools. For recent results in this direction, where asymptotic expansions are performed for $a\to\infty$, see~\cite{SVxx}.
%------------
%-----------

\begin{figure}[h!]
 \begin{minipage}[t]{.5\linewidth}
\centering\epsfig{figure=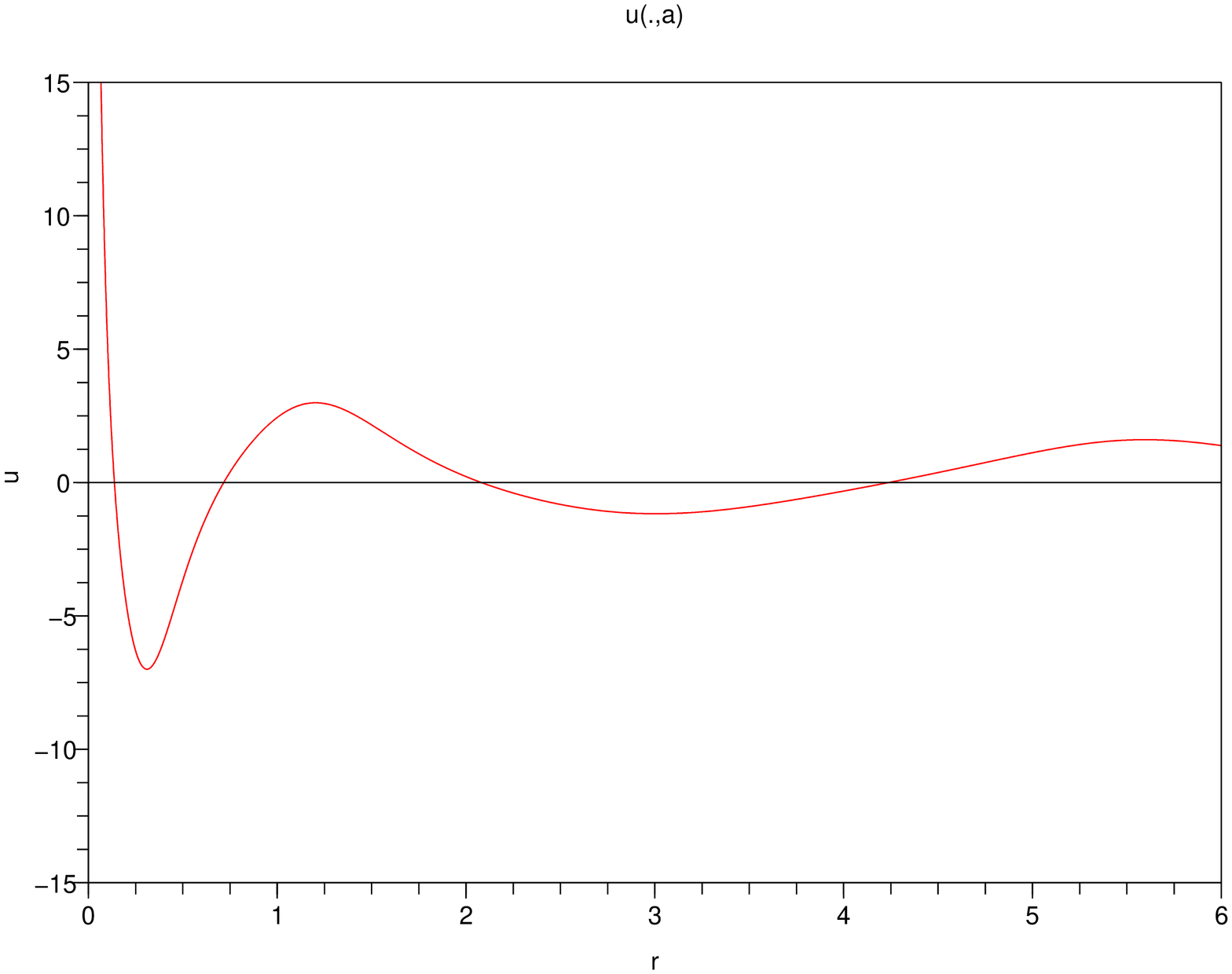,width=4cm,angle=270}
 \end{minipage}\hfill
 \begin{minipage}[t]{.5\linewidth}
\centering\epsfig{figure=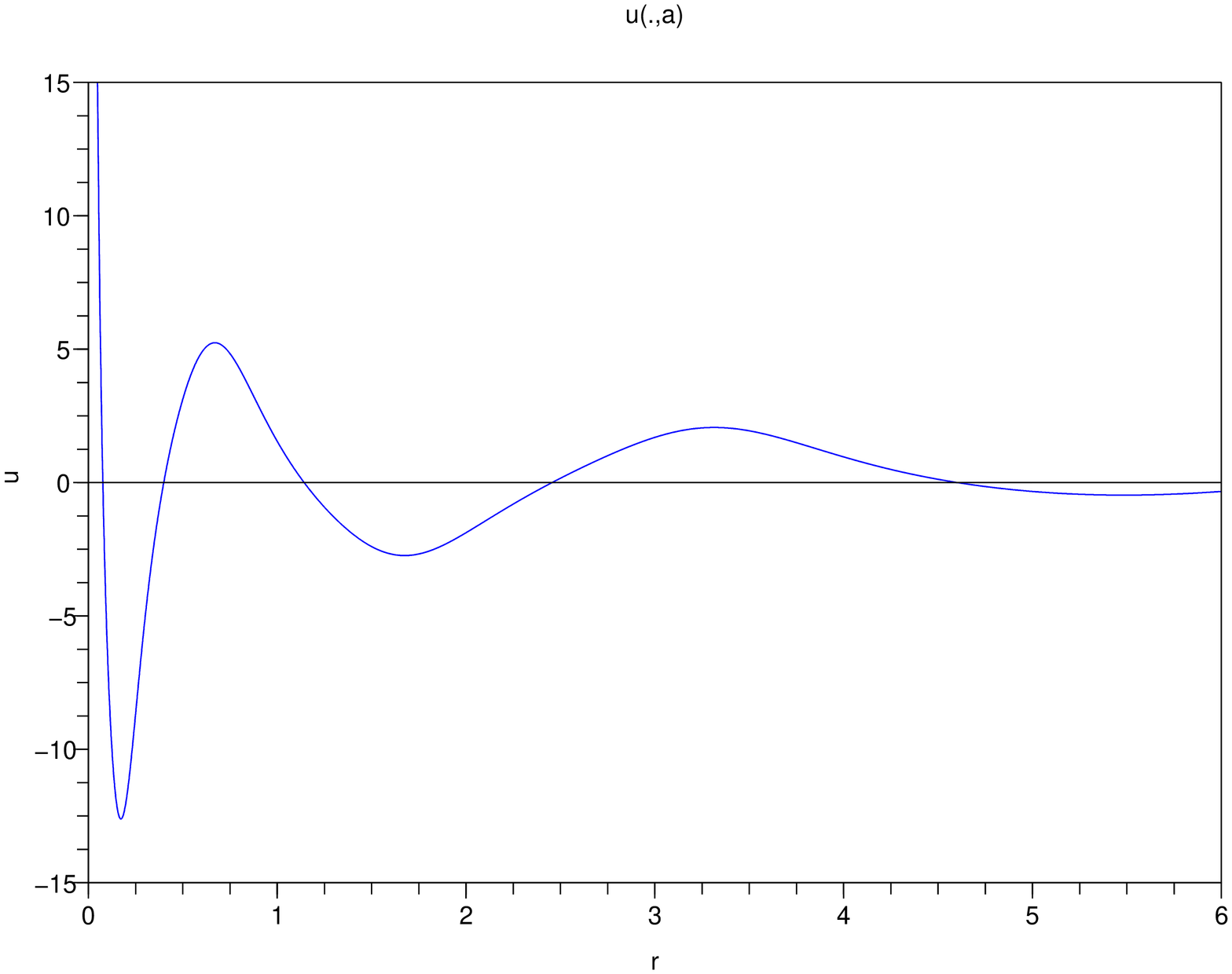,width=4cm,angle=270}
 \end{minipage}
   \caption{Positivity set of $u(.,a)$ with two ($a=50$, left) and three ($a=90$, right) connected components ($d=3$).}\label{fig:prc6}
\end{figure}

\quid{If $a>a_c$ is large enough, $u(a,\cdot)$ may have several zeros, see Figures~\ref{fig:prc6} and each hump corresponds to a solution. It is possible to construct self-similar blowing-up solutions of any mass?}

\quid{The stability of blowing-up solutions is also of interest but seems yet unclear according to numerical simulations performed in \cite{SC08}.}

%------------
%%%%%%%%%%%%%%%%%%%%%%%%%%%%%%%%%%%%%%%%%%%%%%
\section{Concluding remarks}\label{sec:conclusion}
The Patlak-Keller-Segel models have attracted much attention these last years. The literature is vast and drastically increasing. Much work is also currently carried out on the original Patlak-Keller-Segel system with a parabolic equation on the chemo-attractant where the authors prove global-in-time existence for a mass less that $8\pi$~\cite{corriascalvez} but also that there exist global-in-time solutions for larger masses~\cite{bilerdolbeault}. There exist many variants of the presented models with prevention of overcrowding~\cite{MR1826309,schmeiserdolak} or with non-linear chemo-sensitivity~\cite{bertozzigroup,nouvosugiyama}, etc. The Patlak-Keller-Segel model has recently been used as a basis for more complete models~\cite{vincentgroup}. This review is dedicated to the parabolic-elliptic Patlak-Keller-Segel model and tries to describe the progress made through energy and functional inequalities methods in the idea of~\cite{villani,AmbrosioGigliSavare02p,Carrillo-McCann-Villani03}. For more complete reviews see~\cite{H03,H03a,perthame,HP}.

From the author's point of view, the most challenging question is the understanding of the blowup. And in this direction progress is still to be made. We are now at a point where we need to develop new methods to address those questions. The answer could come from interaction with the non-linear Schr\"odinger equation (NLS) and the unstable thin-film equation (UTF). Indeed, the Patlak-Keller-Segel, the NLS and the UTF equation have two levels of criticality. The first level is given by the homogeneity of the ``attractive'' and ``repulsive'' terms in each problem. In our particular case, this refers to the aggregation versus diffusion mechanisms. In NLS it is the balance between dispersion and nonlinear attraction. As seen above, the balance happens precisely for our chosen exponent $m=m_d$. In the NLS equation this happens for the so-called \emph{pseudo-conformal} non-linearity, see~\cite{sulemsulem} or~\cite[Chapter 6]{Ca03}. In the UTF equation this happens in the so-called \emph{marginal} case, see~\cite{BP1,BP2}. In theses three equations, a second level of critically occurs when the attractive and repulsive are balanced. In that particular case, and for the three models, there exists a critical value $M_c$ of the mass which is the maximum value of the mass below which the solutions exist globally in time, see~\cite{Ca03,Weinstein83,MR} for the pseudo-conformal NLS equation and~\cite{BP2,SP} for the marginal UTF equation. Note that mass refers to the total number of particles for the NLS equation and the $\LL^1$-norm for \eqref{eq:spf} and the UTF equation.

Let us also point out that in all these three problems, the virial method is an elegant way to prove that there are solutions which blowup above the critical mass, but it does not give any hint on the mechanism of the blowup. However, in \cite{MR}, for the NLS equation, the result goes further and clarifies the blow-up for super-critical masses close to critical.

The collapsing solutions to the NLS equation have the form of a rescaled ground state soliton and the blowup profile~\eqref{mechablo} has the form of a rescaled steady-state solution of the Keller-Segel system. The scaling of the leading term in the collapsing solutions has the same order $\sqrt{T-t}$ in both NLS and KS.\medskip

%%%%%%%%%%%%%%%%%%%%%%%%%%%%%%%%%%%%%%%%%%%%%%%%%%%%%%%%%%
%%%%%%%%%%%%%%%%%%%%%%%%%%%%%%%%%%%%%%%%%%%%%%%%%%%%%%%%%%
%%%%%%%%%%%%%%%%%%%%%%%%%%%%%%%%%%%%%%%%%%%%%%%%%%%%%%%%%%
%%% Acknowledgement  %%%%%%%%%%%%%%%%%%%%%%%%%%%%%%%%%%%%%
%%%%%%%%%%%%%%%%%%%%%%%%%%%%%%%%%%%%%%%%%%%%%%%%%%%%%%%%%%
\noindent {\bf Acknowledgements.-} The author thanks J. Dolbeault and Ph. Lauren\c cot for their attentive reading of this manuscript and for their suggestions.

%%%%%%%%%%%%%%%%%%%%%%%%%%%%%%%%
%%%%%%%%%%%%%%%%%%%%%%%%%%%%%%%%
%%%%%%%%%%%%%%%%%%%%%%%%%%%%%%%%%%%%%%%%%%%%%%%%%%%%%%%%%%%%%%%
%%%  Copyright   %%%%%%%%%%%%%%%%%%%%%%%%%%%%%%%%%%%%%%%%%%%%%%
%%%%%%%%%%%%%%%%%%%%%%%%%%%%%%%%%%%%%%%%%%%%%%%%%%%%%%%%%%%%%%%
\bigskip\noindent{\small This paper is under the Creative Commons
licence Attribution-NonCommercial-ShareAlike 2.5.}
%%%%%%%%%%%%%%%%%%%%%%%%%%%%%%%%%%%%%%%%%%%%%%%%%%%%%%%%%%%%%%%%%
%%%%%%%%%%%%%%%%%%%%%%%%%%%%%%%%%%%%%%%%%%%%%%%%%%%%%%%%%%%%%%%%%

%%%%%%%%%%%%%%%%%%%%%%%%%%%%%%%%%%%%%%%%%%%%%%%%%%%%%%%%%%
%%% Bibliography %%%%%%%%%%%%%%%%%%%%%%%%%%%%%%%%%%%%%%%%%
%%%%%%%%%%%%%%%%%%%%%%%%%%%%%%%%%%%%%%%%%%%%%%%%%%%%%%%%%%

%%%%%%%%%%%%%%%%%%%%%%%%%%%%%%%%
%%%%%%%%%%%%%%%%%%%%%%%%%%%%%%%%

\begin{thebibliography}{99}

\bibitem{AmbrosioGigliSavare02p}
{\sc L.~A.~Ambrosio, N.~Gigli, and G.~Savar\'e}, {\em Gradient flows
in metric spaces and in the space of probability measures},
Lectures in Mathematics, Birkh\"auser, 2005.

\bibitem{bertozzigroup}
{\sc J. Bedrossian, N. Rodr\'iguez and A. Bertozzi}, {\em Local and global well-posedness for aggregation equations and Patlak-Keller-Segel models with degenerate diffusion}, Nonlinearity, 24 (2011), pp.~1683--1715.

\bibitem{BP1}
{\sc A.~L. Bertozzi and M.~C. Pugh}, {\em Long-wave instabilities
and saturation in thin film equations}, Comm. Pure Appl. Math., 51
(1998), pp.~625--661.

\bibitem{BP2}
\leavevmode\vrule height 2pt depth -1.6pt width 23pt, {\em Finite-time blow-up of
solutions of some long-wave unstable thin film equations}, Indiana
Univ. Math. J., 49 (2000), pp.~1323--1366.

\bibitem{bilerdolbeault}
{\sc P. Biler, L. Corrias and J. Dolbeault}, {\em Large mass self-similar solutions of the parabolic-parabolic Keller-Segel model of chemotaxis}, J. of Math. Biol., 63 (2011), pp.~1--32.

\bibitem{bkln}
{\sc P.~Biler, G.~Karch, P.~Lauren{\c{c}}ot and T.~Nadzieja}, {\em The $8\pi$-problem for radially symmetric solutions of a chemotaxis model in the plane}, Math. Methods Appl. Sci., 29 (2006), pp.~1563--1583.

\bibitem{biler-nadzieja}
{\sc P.~Biler and T.~Nadzieja}, {\em Global and exploding
  solutions in a model of self-gravitating systems}, Rep. Math. Phys., 52
  (2003), pp.~205--225.

\bibitem{BCC}
{\sc A. Blanchet, V.~Calvez and J.A. Carrillo}, {\em Convergence
of the mass-transport steepest descent scheme for the subcritical
{P}atlak-{K}eller-{S}egel model}, SIAM J. Numer. Anal., 46 (2008),
pp.~691--721.

\bibitem{BCCxx}
{\sc A.~Blanchet, E.~Carlen and J.A.~Carrillo}, {\em Functional inequalities, thick tails and asymptotics for the critical mass Patlak-Keller-Segel model},
preprint, arXiv:{1009.0134}. 

\bibitem{BCL09}
{\sc A.~Blanchet, J.A.~Carrillo and Ph.~Lauren\c cot}, {\em Critical mass for a Patlak-Keller-Segel model with degenerate diffusion in higher dimensions}, Calc. Var. Partial Differential Equations \textbf{35} (2009), pp.~133--168.

\bibitem{BCM}
{\sc A.~Blanchet, J.~A.~Carrillo and N.~Masmoudi}, {\em Infinite
Time Aggregation for the Critical Patlak-Keller-Segel model in
$\RR^2$}, Comm. Pure Appl. Math., 61 (2008), pp.~1449--1481.

\bibitem{BDEF}
{\sc A.~Blanchet, J.~Dolbeault, M.~Escobedo and J.~Fernandez}, {\em Asymptotic behaviour for small mass in the two-dimensional parabolic-elliptic Keller-Segel model}, Journal of Mathematical Analysis and Applications, 361 (2010), pp.~533--542..

\bibitem{BDP}
{\sc A.~Blanchet, J.~Dolbeault and B.~Perthame}, {\em
Two-dimensional {K}eller-{S}egel model: optimal critical mass and qualitative properties of the solutions}, Electron. J. Differential Equations 44 (2006), 32 pp. (electronic).

\bibitem{bl}
{\sc A.~Blanchet and Ph.~Lauren\c cot}, {\em Finite mass self-similar blowing-up solutions of a chemotaxis system with non-linear diffusion}, To appear in Communications on Pure and Applied Analysis (2010).

\bibitem{brenner}
{\sc M.~P.~Brenner, L.~S.~Levitov and E.~O.~Budrene}, {\em Physical Mechanisms for Chemotactic Pattern Formation by Bacteria}, Biophysical Journal 74 (1998), pp.~1677--1693.

\bibitem{schmeiserdolak}
{\sc M. Burger, Y. Dolak-Struss, and C. Schmeiser}, {\em Asymptotic analysis of an advection-dominated chemotaxis model in multiple spatial dimensions}, Commun. Math. Sci., 6 (2008) pp.~1--28. 

\bibitem{calvezcarrillo}
{\sc V. Calvez and J.~A. Carrillo}, {\em Volume effects in the Keller-Segel model: energy estimates preventing blow-up}, J. Math. Pure et Appl. 86 (2006), pp.~155--175.

\bibitem{corriascalvez}
{\sc V. Calvez and L. Corrias}, {\em The parabolic-parabolic Keller-Segel model in $\R^2$}, Commun. Math. Sci. 6 (2008), pp.~417--447. 

\bibitem{vincentgroup}
{\sc V. Calvez, R. Hawkins, N. Meunier and R. Voituriez}, {\em Analysis of a non local model for spontaneous cell polarisation}, arXiv:{1105.4429}.

\bibitem{CFxx}
{\sc E. Carlen and A. Figalli}, {\em Stability for a GNS inequality and the Log-HLS inequality,
with application to the critical mass Keller-Segel equation}, preprint, arXiv:{1107.5976}.

\bibitem{Carrillo-McCann-Villani03}
{\sc J.~A. Carrillo, R.~J. McCann and C. Villani}, {\em Kinetic
equilibration rates for granular media and related equations:
entropy dissipation and mass transportation estimates}, Rev.
Matem\'atica Iberoamericana, 19 (2003), pp.~1--48.

\bibitem{Ca03}
{\sc T.~Cazenave}, {\em Semilinear Schr\"odinger
equations}, Courant Lecture Notes in Mathematics, \textbf{10}, New
York University, Courant Institute of Mathematical Sciences, New York;
American Mathematical Society, Providence, RI, 2003.

\bibitem{Chandrasekhar}
{\sc S. Chandrasekhar}, {\em An Introduction to the Study of Stellar Structure}, New York, Dover, 1967.

\bibitem{ChMannella}
{\sc P.-H.~Chavanis and R. Mannella}, {\em Self-gravitating Brownian particles in two dimensions: the case of $N = 2$ particles}, The Eur. Phys. J. B, 78 (2010), pp.~139--165.

\bibitem{ChS04}
{\sc P.-H.~Chavanis and C.~Sire}, {\em Anomalous diffusion and collapse of self-gravitating Langevin particles in $D$ dimensions}, Phys. Rev. E, 69 (2004), 016116.

\bibitem{MR632161}
{\sc S.~Childress and J.~K. Percus}, {\em Nonlinear aspects of chemotaxis},
  Math. Biosci., 56 (1981), pp.~217--237.

\bibitem{CPZ}
{\sc L.~Corrias, B.~Perthame, and H.~Zaag}, {\em Global solutions of
  some chemotaxis and angiogenesis systems in high space dimensions}, Milan J.
  Math., 72 (2004), pp.~1--29.

\bibitem{DD} {\sc M.~Del Pino and J.~Dolbeault}, {\em Best constants
for Gagliardo-Nirenberg inequalities and applications to nonlinear
diffusions}, J. Math. Pures Appl., 81 (2002), pp.~847--875.

\bibitem{DP}
{\sc J.~Dolbeault and B.~Perthame}, {\em
Optimal critical mass in the two-dimensional {K}eller-{S}egel model in {$\mathbb R\sp 2$}}, C. R. Math. Acad. Sci. Paris, 339 (2004), pp.~611--616.

\bibitem{DS}
{\sc J.~Dolbeault and C.~Schmeiser}, {\em The two-dimensional
{K}eller-{S}egel model after blow-up}, Disc. Cont. Dynam. Systems
B, 25 (2009), pp.~109--121.

\bibitem{MR1415081}
{\sc M.~A. Herrero and J.~J.~L. Vel{\'a}zquez}, {\em Singularity patterns in a
  chemotaxis model}, Math. Ann., 306 (1996), pp.~583--623.

\bibitem{MR1826309}
{\sc T.~Hillen and K.~Painter}, {\em Global existence for a parabolic
  chemotaxis model with prevention of overcrowding}, Adv. in Appl. Math., 26
  (2001), pp.~280--301.

\bibitem{HP}
{\sc T. Hillen and K. Painter}, {\em A user's guide to PDE
models for chemotaxis}, Journal of Mathematical Biology, 58
(2009), pp.~183--217.

\bibitem{H03a}
{\sc D.~Horstmann}, {\em From 1970 until
  present: the {K}eller-{S}egel model in chemotaxis and its consequences. {I}},
  Jahresber. Deutsch. Math.-Verein., 105 (2003), pp.~103--165.

\bibitem{H03}
\leavevmode\vrule height 2pt depth -1.6pt width 23pt, {\em From 1970 until
  present: the {K}eller-{S}egel model in chemotaxis and its consequences.
  {II}}, Jahresber. Deutsch. Math.-Verein., 106 (2004), pp.~51--69.

\bibitem{JL92}
{\sc W.~J\"{a}ger and S.~Luckhaus}, {\em On explosions of
solutions to a system of partial differential equations modelling
chemotaxis}, Trans. Amer. Math. Soc., 329 (1992), pp.~819--824.

\bibitem{souplet}
{\sc N. Kavallaris and P. Souplet}, {\em Grow-up rate and refined asymptotics for a two-dimensional Patlak-Keller-Segel model in a disk}, SIAM J. Math. Anal., 40 (2008/09), pp.~1852--1881.

\bibitem{Keller-Segel-70}
{\sc E.~F. Keller and L.~A. Segel}, {\em Initiation of slide mold
aggregation viewed as an instability}, J. Theor. Biol., 26 (1970), pp.~399--415.

\bibitem{Kowalczyk04}
{\sc R.~Kowalczyk}, {\em Preventing blow-up in a chemotaxis
model}, J. Math. Anal. Appl. 305 (2005), pp.~566--588.

\bibitem{lm}
{\sc C.~Lederman and P.~A.~Markowich}, {\em On fast-diffusion
equations with infinite equilibrium entropy and finite equilibrium
mass}, Comm. Partial Differential Equations, 28 (2003),
pp.~301--332.

\bibitem{nouvosugiyama}
{\sc S. Luckhaus, Y. Sugiyama, J. J. L. Vel\'azquez}, {\em Measure valued solutions of the 2D Keller-Segel system}, arXiv:1011.0282.

\bibitem{Luxx} 
{\sc P.~M.~Lushnikov}, {\em Critical chemotactic collapse}, Phys. Lett.~A \textbf{374} (2010), pp.~1678--1685. 

\bibitem{medina}
{\sc M. A. Herrero, E. Medina and J. J. L. Vel\'azquez}, {\em Self-similar blow-up for a reaction-diffusion system}, J. of Comp. and Appl. Math., 97 (1998), pp.~99--119. 

\bibitem{MR}
{\sc F. Merle and P. Rapha\"el}, {\em On universality of
blow-up profile for {$L\sp 2$} critical nonlinear {S}chr\"odinger
equation}, Invent. Math., 156 (2004), pp.~565--672.

\bibitem{nagai}
{\sc T. Nagai}, {\em Blow-up of radially symmetric
solutions to a chemotaxis system}, Adv. Math. Sci. Appl. 5 (1995), pp. 581--601.

\bibitem{nagaisenbasuzuki}
{\sc T. Nagai, T. Senba, T. Suzuki}, {\em Chemotactic collapse in a parabolic system of mathematical biology}, Hiroshima Math. J, 30 (2000), pp.~463--497.

\bibitem{nanjundiah}
{\sc V.~Nanjundiah}, {\em Chemotaxis, signal relaying and aggregation
  morphology}, Journal of Theoretical Biology, 42 (1973), pp.~63--105.

\bibitem{MR1052381}
{\sc T.~Padmanabhan}, {\em Statistical mechanics of gravitating systems}, Phys.
  Rep., 188 (1990), pp.~285--362.

\bibitem{Pat53}
{\sc C.~S. Patlak}, {\em Random walk with persistence and external
bias}, Bull. Math. Biophys., 15 (1953), pp.~311--338.

\bibitem{perthame}
{\sc B.~Perthame}, {\em Transport equations in biology}, Frontiers in Mathematics, Birkh\"auser Verlag, Basel, 2007.

\bibitem{SP}
{\sc D. Slep{\v{c}}ev and M.~C. Pugh}, {\it Selfsimilar blowup of
unstable thin-film equations}, Indiana Univ. Math. J., 54 (2005),
pp.~1697--1738.

\bibitem{SC02} 
{\sc C.~Sire and P.-H.~Chavanis}, {\em Thermodynamics and collapse of self-gravitating Brownian particles in D dimensions}, Phys. Rev. E, 66 (2002), 046133.

\bibitem{SC08} {\sc C.~Sire and P.-H.~Chavanis}, {\em Critical dynamics of self-gravitating Langevin particles and bacterial populations},  Phys. Rev. E, 78 (2008), 061111.

\bibitem{S1}
{\sc Y. Sugiyama}, {\em Global existence in sub-critical cases and
finite time blow-up in super-critical cases to degenerate
Keller-Segel systems}, Differential Integral Equations 19 (2006),
pp.~841--876.

\bibitem{S2}
{\sc Y. Sugiyama}, {\em Application of the best constant of the
Sobolev inequality to degenerate Keller-Segel models}, Adv.
Differential Equations 12 (2007), pp.~121--144.

\bibitem{SVxx}
{\sc Y. Sugiyama and J.~J.~L. Vel{\'a}zquez}, {\em Self-similar blow-up with a continuous range of values of the aggregated mass for a degenerate Keller-Segel system},
Advances in Differential Equations, 16 (2011), pp.~85--112. 

\bibitem{sulemsulem}
{\sc C. Sulem and P.~L. Sulem}, {\em The nonlinear {S}chr\"odinger
equation}, Applied Mathematical Sciences 139, Springer-Verlag, New
York, 1999.

\bibitem{ST09}
{\sc T. Suzuki and R.~Takahashi}, {\em Degenerate parabolic equation with critical exponent derived from the kinetic theory, I, generation of the weak solution}, Adv. Differential Equations, 14 (2009), pp.~433--476.

\bibitem{velazquez}
{\sc J.~J.~L. Vel{\'a}zquez}, {\em Stability of some mechanisms of chemotactic
  aggregation}, SIAM J. Appl. Math., 62 (2002), pp.~1581--1633 (electronic).

\bibitem{MR2068667}
\leavevmode\vrule height 2pt depth -1.6pt width 23pt, {\em Point dynamics in a
  singular limit of the {K}eller-{S}egel model. {I}. {M}otion of the
  concentration regions}, SIAM J. Appl. Math., 64 (2004), pp.~1198--1223
  (electronic).

\bibitem{MR2068668}
\leavevmode\vrule height 2pt depth -1.6pt width 23pt, {\em Point dynamics in a
  singular limit of the {K}eller-{S}egel model. {II}. {F}ormation of the
  concentration regions}, SIAM J. Appl. Math., 64 (2004), pp.~1224--1248
  (electronic).

\bibitem{villani}
{\sc C. Villani},
{\em Topics in optimal transportation},
Graduate Studies in Mathematics Vol. 58, Amer. Math. Soc, Providence, 2003.

\bibitem{Weinstein83}
{\sc M.~I. Weinstein}, {\em Nonlinear Schr\"odinger equations and
sharp interpolation estimates}, Comm. Math. Phys. 87 (1983),
pp.~567--576.
\end{thebibliography}
\end{document}